\documentclass[12pt]{amsart}

\usepackage{amsmath}
\usepackage{amscd}
\usepackage{amssymb}
\usepackage{latexsym}
\usepackage{color}
\usepackage{mathrsfs}
\usepackage{hyperref}
\usepackage[all]{xy}
\usepackage{tikz-cd}
\usepackage{stmaryrd}

\newtheorem{theorem}{Theorem}[section] 
\newtheorem*{theorem*}{Theorem}
\newtheorem{lemma}[theorem]{Lemma}
\newtheorem*{lemma*}{Lemma}
\newtheorem{corollary}[theorem]{Corollary}
\newtheorem*{corollary*}{Corollary}
\newtheorem{proposition}[theorem]{Proposition}
\newtheorem*{proposition*}{Proposition}

\newtheorem{remark}[theorem]{Remark}
\newtheorem{question}[theorem]{Question}
\newtheorem{definition}[theorem]{Definition}

\newtheorem{example}[theorem]{Example}
\newcommand{\bgl}{\begin{equation}} 
\newcommand{\egl}{\end{equation}}
\newcommand{\bgloz}{\begin{equation*}} 
\newcommand{\egloz}{\end{equation*}}
\newcommand{\bgln}{\begin{eqnarray}} 
\newcommand{\egln}{\end{eqnarray}}
\newcommand{\bglnoz}{\begin{eqnarray*}} 
\newcommand{\eglnoz}{\end{eqnarray*}}
\newcommand{\btheo}{\begin{theorem}}
\newcommand{\etheo}{\end{theorem}}
\newcommand{\btheooz}{\begin{theorem*}}
\newcommand{\etheooz}{\end{theorem*}}
\newcommand{\blemma}{\begin{lemma}}
\newcommand{\elemma}{\end{lemma}}
\newcommand{\blemmaoz}{\begin{lemma*}}
\newcommand{\elemmaoz}{\end{lemma*}}
\newcommand{\bproof}{\begin{proof}}
\newcommand{\eproof}{\end{proof}}
\newcommand{\bbew}{\begin{beweis}}
\newcommand{\ebew}{\end{beweis}}
\newcommand{\bremark}{\begin{remark}\em}
\newcommand{\eremark}{\end{remark}}
\newcommand{\bquestion}{\begin{question}\em}
\newcommand{\equestion}{\end{question}}
\newcommand{\bdefin}{\begin{definition}}
\newcommand{\edefin}{\end{definition}}
\newcommand{\bprop}{\begin{proposition}}
\newcommand{\eprop}{\end{proposition}}
\newcommand{\bpropoz}{\begin{proposition*}}
\newcommand{\epropoz}{\end{proposition*}}
\newcommand{\bcor}{\begin{corollary}}
\newcommand{\ecor}{\end{corollary}}
\newcommand{\bcoroz}{\begin{corollary*}}
\newcommand{\ecoroz}{\end{corollary*}}
\newcommand{\bfa}{\begin{cases}} 
\newcommand{\efa}{\end{cases}}
\newcommand{\bexample}{\begin{example}\em}
\newcommand{\eexample}{\end{example}}
\newcommand{\cB}{\mathcal B}

\newcommand{\cG}{\mathcal G}

\newcommand{\cM}{\mathcal M}

\newcommand{\cO}{\mathcal O}
\newcommand{\cP}{\mathcal P}

\newcommand{\cS}{\mathcal S}

\def\Cz{\mathbb{C}}

\def\Nz{\mathbb{N}}

\def\Tz{\mathbb{T}}
\def\Zz{\mathbb{Z}}

\def\1z{\mathbb{1}}

\newcommand{\fX}{\mathfrak X}
\newcommand{\bfg}{{\bf g}}

\newcommand{\bfr}{{\bf r}}

\newcommand{\scH}{\mathscr H}
\newcommand{\an}[1]{``#1''} 
\newcommand{\ti}{\tilde}

\newcommand{\lori}{\longrightarrow}

\newcommand{\ma}{\mapsto} 
\newcommand{\mafr}{\mapsfrom} 
\newcommand\onto{\twoheadrightarrow} 
\newcommand\into{\hookrightarrow} 

\newcommand{\ve}{\varepsilon}

\def\SEMI{\mbox{$\times\kern-2pt\vrule height5pt width.6pt \kern3pt $}}

\newcommand{\Hom}{{\rm Hom}\,}

\newcommand{\Aut}{{\rm Aut}\,}

\newcommand{\Spec}{{\rm Spec\,}} 

\newcommand{\id}{{\rm id}}

\newcommand{\Ad}{{\rm Ad\,}}

\renewcommand{\ker}{{\rm ker}\,}

\newcommand{\abs}[1]{\left|#1\right|} 
\newcommand{\norm}[1]{\left\|#1\right\|} 
\newcommand{\defeq}{\mathrel{:=}} 

\newcommand{\dop}{\text{: }} 
\newcommand{\ilim}{\varinjlim} 
\newcommand{\plim}{\varprojlim} 

\newcommand{\supp}{{\rm supp}\,}

\newcommand{\dom}{{\rm dom}\,}
\newcommand{\ran}{{\rm ran}\,}

\newcommand{\ue}{{\approx_{\mathrm{u}}}}

\newcommand{\lge}{\left\{} 
\newcommand{\rge}{\right\}} 
\newcommand{\lru}{\left(} 
\newcommand{\rru}{\right)} 
\newcommand{\rukl}[1]{\lru #1 \rru} 
\newcommand{\gekl}[1]{\lge #1 \rge} 
\newcommand{\vp}{\varphi}
\newcommand{\Rho}{{\rm P}}
\newcommand{\menge}[2]{\gekl{ #1 \dop #2 }} 
\newcommand{\Ext}{{\rm Ext}\,}
\newcommand{\bG}{{\bar G}}
\newcommand{\bSigma}{{\bar \Sigma}}
\newcommand{\bX}{{\bar X}}
\newcommand{\etale}{{\'e}tale}
\title{Cartan subalgebras and the UCT problem, II}

\author{Sel\c{c}uk Barlak}
\address{Aschaffenburg\\Germany}
\email{selcuk@barlak.de}

\author{Xin Li}
\address{School of Mathematics and Statistics\\
University of Glasgow\\
University Place\\
Glasgow G12 8QQ\\
United Kingdom}
\email{Xin.Li@glasgow.ac.uk}

\subjclass[2010]{Primary 46L05, 46L40; Secondary 46L80, 19K35}

\thanks{The first named author is supported by  the Villum Fonden project grant `Local and global
structures of groups and their algebras' (2014–-2018).}
\thanks{The second named author is supported by EPSRC grant EP/M009718/1.}

\begin{document}

\begin{abstract}
We study the connection between the UCT problem and Cartan subalgebras in C*-algebras. The UCT problem asks whether every separable nuclear C*-algebra satisfies the UCT, i.e., a noncommutative analogue of the classical universal coefficient theorem from algebraic topology. This UCT problem is one of the remaining major open questions in the structure and classification theory of simple nuclear C*-algebras. Since the class of separable nuclear C*-algebras is closed under crossed products by finite groups, it is a natural and important task to understand the behaviour of the UCT under such crossed products. We make a contribution towards a better understanding by showing that for certain approximately inner actions of finite cyclic groups on UCT Kirchberg algebras, the crossed products satisfy the UCT if and only if we can find Cartan subalgebras which are invariant under the actions of our finite cyclic groups.

We also show that the class of actions we are able to treat is big enough to characterize the UCT problem, in the sense that every such action (even on a particular Kirchberg algebra, namely the Cuntz algebra $\cO_2$) leads to a crossed product satisfying the UCT if and only if every separable nuclear C*-algebra satisfies the UCT.

Our results rely on a new construction of Cartan subalgebras in certain inductive limit C*-algebras. This new tool turns out to be of independent interest. For instance, among other things, the second author has used it to construct Cartan subalgebras in all classifiable unital stably finite C*-algebras.
\end{abstract}

\maketitle



\section{Introduction}

\setlength{\parindent}{0cm} \setlength{\parskip}{0cm}

This paper studies the connection between the UCT problem -- one of the remaining major open questions in the structure and classification theory of simple nuclear C*-algebras -- and Cartan subalgebras of C*-algebras -- a concept of recent interest which builds bridges between C*-algebras, topological dynamics, and geometric group theory.
\setlength{\parindent}{0cm} \setlength{\parskip}{0.5cm}

The universal coefficient theorem (UCT) for C*-algebras was introduced by Rosenberg and Schochet in \cite{RSch}. The idea was to formulate a noncommutative analogue of the classical universal coefficient theorem in algebraic topology, very much in line with the general philosophy of viewing C*-algebra theory as noncommutative topology. It facilitates the computation of Kasparov's bivariant $K$-theory \cite{Kas} in terms of topological $K$-theory for C*-algebras. More precisely, a separable C*-algebra $A$ is said to satisfy the UCT if for every separable C*-algebra $A'$ the sequence
$$
 0 \to \Ext(K_*(A),K_{*-1}(A')) \to KK_*(A,A') \to\Hom(K_*(A),K_*(A')) \to 0
$$
is exact, where the right hand map is the natural one and the left hand map is the inverse of a map that is always defined. Equivalently, a separable C*-algebra satisfies the UCT if and only if it is $KK$-equivalent to a commutative C*-algebra. The class of separable nuclear C*-algebras satisfying the UCT coincides with the smallest class of separable nuclear C*-algebras that contains the algebra of complex numbers and is closed under countable inductive limits, the two out of three property for extensions, and $KK$-equivalences, see \cite{Bla}.
\setlength{\parindent}{0.5cm} \setlength{\parskip}{0cm}

Although there exist exact non-nuclear C*-algebras that do not satisfy the UCT, \cite{S}, it is still open whether all separable nuclear C*-algebras satisfy the UCT. This question is often referred to as the UCT problem. Due to the recent breakthrough results in the classification program of separable simple nuclear C*-algebras, \cite{EGLN,GLN,TWW}, the UCT problem is receiving renewed interest. In fact, the necessary UCT assumption for classifiable C*-algebras remains mysterious. A positive solution to the UCT problem would show that the UCT assumption is redundant and thereby completely clarify the situation.

Since the class of separable nuclear C*-algebras is closed under crossed products by countable (discrete) amenable groups, it is a natural and important question whether the UCT is preserved under such crossed products. We make a contribution towards a better understanding by recasting this question in terms of Cartan subalgebras in C*-algebras.
\setlength{\parindent}{0cm} \setlength{\parskip}{0.5cm}

The concept of Cartan subalgebras in C*-algebras was introduced and first studied in \cite{Ku1,R}, and builds a bridge between C*-algebras, topological dynamics, and geometric group theory. On the one hand, Cartan subalgebras provide geometric descriptions for C*-algebras and hence introduce ideas from topology and dynamical systems to the study of C*-algebras. On the other hand, Cartan subalgebras and the closely related notion of continuous orbit equivalence allow us to apply methods from functional analysis and operator algebras in topological dynamics and geometric group theory (see for instance \cite{Li16,Li17,Li18}). 
\setlength{\parindent}{0.5cm} \setlength{\parskip}{0cm}

These developments led to a renewed interest in Cartan subalgebras, and this interest was enhanced by the observation that Cartan subalgebras are closely related to the UCT problem (for instance, the recent Oberwolfach Mini-Workshop \cite{OWR17} was dedicated to this theme). The starting point was Tu's striking result \cite{Tu} that amenable groupoids provide a rich source for nuclear C*-algebras in the UCT class. Building on Tu's work, the authors have proved that this also holds for twisted amenable, \etale{} groupoids, see \cite{BL}. 
Therefore, employing Renault's remarkable characterisation of separable C*-algebras with Cartan subalgebras as twisted \etale{} groupoids \cite{R}, it follows that all separable nuclear C*-algebras with Cartan subalgebras satisfy the UCT.
At the same time, by work of Spielberg~\cite{Sp07a} (or results of Katsura~\cite{Kat08c} and Yeend \cite{Y1,Y2}), UCT Kirchberg algebras are known to admit a Cartan subalgebra. Combining these results with the Kirchberg-Phillips classification theorem \cite{Kir, Phi}, the UCT problem turns out to be equivalent to the question whether every Kirchberg algebra has a Cartan subalgebra. 

This prompts the natural question whether an analogous characterization of the UCT problem is possible using stably finite C*-algebras. In this context, we would like to mention the recent work of Deeley, Putnam and Strung, \cite{DPS}, on a realisation of the Jiang-Su algebra \cite{JS,RW} via an \etale{} equivalence relation. Moreover, building on Renault's groupoid C*-algebra models for AF algebras \cite{Rbook}, Putnam \cite{Put} recently showed that a large class of classifiable stably finite C*-algebras can be constructed from \etale{} groupoids. Recently, the second author has applied the machinery developed in our paper (Theorem~\ref{indlim--GPD}, see explanations below) to show -- among other things -- that every classifiable unital stably finite C*-algebras has a Cartan subalgebra \cite{Li18'}. Together with the reformulation of the UCT problem in \cite{Dad}, this indeed shows that the UCT problem is equivalent to the question whether every unital separable simple stably finite C*-algebra with finite nuclear dimension has a Cartan subalgebra. In a more general setting, existence and uniqueness of Cartan subalgebras has also been studied in \cite{LR}.
\setlength{\parindent}{0cm} \setlength{\parskip}{0.5cm}

Given that all known classifiable C*-algebras arise from (twisted) groupoids, it is natural to investigate whether dynamics involving such C*-algebras are geometric in nature as well. More concretely, one can ask which discrete group actions on classifiable C*-algebras are (cocycle) conjugate to actions that are induced by twisted groupoid automorphisms. With regards to Renault's characterisation of Cartan subalgebras, our starting point was the following question. 
\bquestion \label{ques:intro}
Let $A$ be a C*-algebra that is classifiable in the sense of the Elliott program and $\Gamma$ a countable discrete amenable group. Which $\Gamma$-actions of $A$ by *-automorphisms preserve some Cartan subalgebra of $A$ globally?
\equestion
\setlength{\parindent}{0cm} \setlength{\parskip}{0cm}

$A$ being nuclear and $\Gamma$ being amenable, the crossed product $A \rtimes_r \Gamma$ will be nuclear, so that a positive solution to the UCT problem would imply that $A \rtimes_r \Gamma$ satisfies the UCT; and the presence of a $\Gamma$-invariant Cartan subalgebra in $A$ would give a conceptual explanation why the crossed product satisfies the UCT because of \cite[Proposition~1.6]{BL}.
\setlength{\parindent}{0cm} \setlength{\parskip}{0.5cm}

In this paper, we restrict attention to actions of finite groups on UCT Kirchberg algebras that either have the Rokhlin property or are approximately representable in the sense of Izumi \cite{Iz} and Nawata \cite{Na}. These classes of actions are dual to each other (see \cite{Iz,Na} and also \cite{BS2,Gar}), so that, by employing Takai duality \cite{Tak}, their study often goes hand in hand. In particular, the rigidity of the Rokhlin property makes both, actions with the Rokhlin property and approximately representable actions, natural classes of actions which are particularly accessible to classification, see for instance \cite{Iz,Iz2}. Moreover, it turns out that it is enough to understand such actions of finite cyclic groups in order to give an equivalent characterization of the UCT problem (see Theorem~\ref{thm:intro UCT}).
\setlength{\parindent}{0.5cm} \setlength{\parskip}{0cm}

As shown in \cite{BL}, a necessary condition for a positive answer to Question~\ref{ques:intro} is that the associated crossed products satisfy the UCT. Using Izumi's classification results in \cite{Iz}, it is shown in the same paper that this condition is also sufficient for outer actions of $\Zz_2$ on the Cuntz algebra $\cO_2$ that are approximately representable. Using a new approach and different techniques explained later in the introduction, we establish a vast generalization of the results in \cite{BL}. In the following, let us formulate our main results.
\btheo \label{thm:intro approx repr}
Let $n \geq 2$. Let $A$ be a unital UCT Kirchberg algebra and $\alpha: \Zz_n \curvearrowright A$ an outer approximately representable action. Assume that $A \rtimes_\alpha \Zz_n$ absorbs the UHF algebra $M_{n^\infty}$ tensorially. Then the following are equivalent:
\begin{itemize}
\item[(i)] $A \rtimes_\alpha \Zz_n$ satisfies the UCT;
\item[(ii)] there exists an inverse semigroup $\cS \subset A$ of partial isometries inducing an isomorphism $A \cong C^*_{tight}(\cS)$ such that 
\begin{itemize}
\item[(1)] $\cS$ is $\alpha$-homogeneous;
\item[(2)] $C^*(E)$ is a Cartan subalgebra with spectrum homeomorphic to the Cantor space in both $A^\alpha$ and $A$;
\item[(3)] unitaries in $A^\alpha$ witnessing approximate representability can be chosen in $C^*(E)$;
\end{itemize}
\item[(iii)] there exists a Cartan subalgebra $C \subset A$ such that $\alpha(C) = C$.
\end{itemize}
\etheo
Here, $C^*_{tight}(\cS)$ is the tight C*-algebra of $\cS$ in the sense of \cite{Ex}. Moreover, we denote by $E$ the idempotent semilattice of $\cS$ and by $C^*(E)$ the sub-C*-algebra of $A$ generated by $E$. The condition in Theorem~\ref{thm:intro approx repr}(ii) can be considered as a quasi-freeness type condition with respect to the inverse semigroup $\cS$. Indeed, if $\cS$ is generated by a Cuntz family of isometries or, more generally, by a Cuntz-Krieger family of partial isometries, then $\alpha$ is quasi-free in the sense of \cite{Ev} and \cite{Z}.

On our way to Theorem~\ref{thm:intro approx repr}, we prove the following structure result for actions with the Rokhlin property on unital UCT Kirchberg algebras.
\btheo \label{thm:intro structure Rokhlin actions}
Let $\Gamma$ be a finite group such that every Sylow subgroup is cyclic. Let $A$ be a unital UCT Kirchberg algebra with $A \cong A \otimes M_{|\Gamma|^\infty}$. Let $\alpha:\Gamma \curvearrowright A$ be an action with the Rokhlin property. 

Then there exists an $\alpha$-invariant inverse semigroup $\cS \subset A$ of partial isometries inducing an isomorphism $A \cong C^*_{tight}(\cS)$ such that $C^*(E)$ is a Cartan subalgebra with totally disconnected spectrum and $\alpha_\gamma(C^*(E)) = C^*(E)$ for all $\gamma \in \Gamma$.
\etheo
In the particular case that $A$ is the Cuntz algebra $\cO_2$ and $\Gamma$ is a cyclic group of prime power order, Theorem~\ref{thm:intro approx repr} applies to all outer strongly approximately inner actions in the sense of Izumi \cite{Iz}, see Corollary~\ref{cor: charac UCT O_2}. We use this to give a new characterisation of the UCT problem. This completes earlier work by the authors in \cite{BL} where a similar characterisation was given for the smaller class of all separable, nuclear C*-algebras that are $KK$-equivalent to their $M_{2^\infty}$-stabilisation. Our result reads as follows.

\btheo \label{thm:intro UCT}
The following statements are equivalent:
\begin{enumerate}
\item[(i)] Every separable, nuclear C*-algebra $A$ satisfies the UCT;
\item[(ii)] for every prime number $p \geq 2$ and every outer strongly approximately inner action $\alpha:\Zz_p \curvearrowright \cO_2$ there exists an inverse semigroup $\cS \subset \cO_2$ of partial isometries inducing an isomorphism $\cO_2 \cong C^*_{tight}(\cS)$ such that 
\begin{itemize}
\item[(1)] $\cS$ is $\alpha$-homogeneous;
\item[(2)] $C^*(E)$ is a Cartan subalgebra with spectrum homeomorphic to the Cantor space in both $\cO_2^\alpha$ and $\cO_2$;
\item[(3)] unitaries in $\cO_2^\alpha$ witnessing approximate representability can be chosen in $C^*(E)$;
\end{itemize}
\item[(iii)] every outer strongly approximately inner $\Zz_p$-action on $\cO_2$ with $p=2$ or $p = 3$ fixes some Cartan subalgebra $B \subseteq \cO_2$ globally.
\end{enumerate}
\etheo
In spirit, Theorem~\ref{thm:intro UCT} is similar to \cite[Theorem~4.17]{BS}, in the sense that both results reformulate the UCT problem in terms of certain $\Zz_p$-actions on $\cO_2$. While the reformulation in \cite[Theorem~4.17]{BS} still refers to the UCT for the crossed product algebras of the $\Zz_p$-actions, the main novelty of Theorem~\ref{thm:intro UCT} is that we can now avoid mentioning the UCT in (ii) and (iii) by discussing existence of Cartan subalgebras which are compatible with the $\Zz_p$-actions and which have further special properties.

The proofs of our main results rely on other results obtained in this paper, which are of independent interest. In particular, we would like to highlight Theorem~\ref{indlim--GPD}. Here we give a new way of constructing Cartan subalgebras in inductive limits of Cartan pairs, under the assumption that each connecting map is realised as the composition of two *-homomorphisms induced by twisted groupoid homomorphisms, where the first one is proper and surjective and the second one is open and injective. This is the first general result on Cartan subalgebras in inductive limit C*-algebras, and turns out to have applications going far beyond the scope of the present paper. For instance, among other things, the second author has used this new tool to construct Cartan subalgebras in all classifiable unital stably finite C*-algebras in \cite{Li18'}.

The second ingredient is Theorem~\ref{Cartan fixed general} which provides a criterion under which a given finite group action with the Rokhlin property on a separable C*-algebra fixes some Cartan subalgebra globally. Here we make use of the concrete model actions described in \cite{BS}. The crucial feature is that it suffices to pass from the Rokhlin action to any family of automorphisms indexed by the group inducing the same map into the automorphism group of the C*-algebra up to approximate unitary equivalence, and find some Cartan subalgebra in the C*-algebra that is fixed globally by all of these automorphisms. Using a classification result by Izumi \cite{Iz2}, it turns out that for UCT Kirchberg algebras it is even enough for this family of automorphisms to induce the same action on $K$-theory.

Lemma~\ref{Lem:Sp+Kat} and Corollary~\ref{Cor:Sp+Kat-stable} together with Lemma~\ref{lifts UCT Kirchberg} now complete the proofs of Theorems~\ref{thm:intro approx repr} and \ref{thm:intro structure Rokhlin actions}. Using Katsura's homological algebra machinery \cite{Kat08a} and Spielberg's graph based models of UCT Kirchberg algebras \cite{Sp07a} and construction of actions on these models \cite{Sp07b},
we show that for every action of a finite group with cyclic Sylow subgroups on a pair of countable abelian groups there exists an action of an \etale{} locally compact groupoid with totally disconnected unit space such that the associated reduced groupoid C*-algebra is a UCT Kirchberg algebra and the induced action on it realises the prescribed action on the level of $K$-theory. 
From this, we conclude the existence of the desired inverse semigroup $\cS$ in Theorems~\ref{thm:intro approx repr} and \ref{thm:intro structure Rokhlin actions}.
\setlength{\parindent}{0cm} \setlength{\parskip}{0.5cm}

The paper is organised as follows. In Section~\ref{Prelim}, we recall the definition of Cartan subalgebras, Renault's characterisation of Cartan pairs in terms of twisted groupoid C*-algebras, as well as the duality result for finite abelian groups for actions with the Rokhlin property and approximately representable actions. The main result of Section~\ref{Sec:indlim} is Theorem~\ref{indlim--GPD}, where Cartan subalgebras in certain inductive limits of Cartan pairs are constructed. In Section~\ref{Sec:Inv}, we prove our main results Theorems~\ref{thm:intro approx repr}, \ref{thm:intro structure Rokhlin actions} and \ref{thm:intro UCT}.

This work was initiated during a visit of the first author at Queen Mary University of London in autumn 2016. He would like to thank the second author and the Geometry and Analysis Group for their kind hospitality during this stay.
\setlength{\parindent}{0.5cm} \setlength{\parskip}{0cm}

\section{Preliminaries}
\label{Prelim}

In the following, all our groupoids are locally compact Hausdorff second countable topological groupoids.

\subsection{From twisted groupoids to C*-algebras}
\label{ss:twistGPD,CSTAR}

Let us recall some definitions from \cite{Ku1,R}. A twisted groupoid is a locally compact groupoid $\Sigma$ together with a free action $\Tz \curvearrowright \Sigma$ such that $\Tz \times \Sigma \to \Sigma, \, (z,\sigma) \ma z \sigma$ is continuous and the quotient $\Sigma / \Tz$ is Hausdorff. In addition, the canonical projection $\Sigma \onto \Sigma / \Tz$ should be a locally trivial principal $\Tz$-bundle. Furthermore, we require that for all $z_1, z_2 \in \Tz$ and $\sigma_1, \sigma_2 \in \Sigma$, $z_1 \sigma_1$ and $z_2 \sigma_2$ are composable if and only if $\sigma_1$ and $\sigma_2$ are composable, and if that is the case, we want to have $(z_1 \sigma_1) (z_2 \sigma_2) = (z_1 z_2) (\sigma_1 \sigma_2)$.

It follows from these axioms that the groupoid structure on $\Sigma$ induces the structure of a locally compact groupoid on $G \defeq \Sigma / \Tz$. Hence we obtain a central extension
$$
  \Tz \times G^{(0)} \overset{\iota}{\rightarrowtail} \Sigma \overset{\pi}{\twoheadrightarrow} G.
$$
Usually, we denote the twisted groupoid as described above by $(G,\Sigma)$. If $G$ is \etale{} (topologically principal), we say that $(G,\Sigma)$ is a twisted \etale{} (topologically principal) groupoid. The quotient map $\pi: \: \Sigma \twoheadrightarrow G$ is often denoted by $\sigma \ma \dot{\sigma}$. Note that $\pi$ is open and closed, as it is the quotient map corresponding to a continuous action of a compact group. It follows that $\pi$ is proper as it is perfect, i.e., $\pi^{-1}(\gekl{\gamma})$ is compact for all $\gamma \in G$.
\setlength{\parindent}{0cm} \setlength{\parskip}{0.5cm}

To construct the reduced C*-algebra of $(G,\Sigma)$ for $G$ \etale{}, let 
$$C(G,\Sigma) \defeq \menge{f \in C(\Sigma,\Cz)}{f(z \sigma) = f(\sigma) \overline{z}},$$
and for $f \in C(G,\Sigma)$, define $\supp(f) \defeq \overline{\menge{\dot{\sigma} \in G}{f(\sigma) \neq 0}}$. Set
$$C_c(G,\Sigma) \defeq \menge{f \in C(G,\Sigma)}{\supp(f) \ {\rm is} \ {\rm compact}}.$$
$C_c(G,\Sigma)$ becomes a *-algebra under the operations
$$
  f^*(\sigma) \defeq \overline{f(\sigma^{-1})}, \ \ \ (f * g)(\sigma) \defeq \sum_{\dot{\tau} \in G_{s(\sigma)}} f(\sigma \tau^{-1}) g(\tau).
$$
Set $X \defeq G^{(0)}$. For $x \in X$, define the Hilbert space
$$\scH_x \defeq \menge{\xi: \: \Sigma_x \to \Cz}{\xi(z \sigma) = \xi(\sigma) \overline{z}, \, \sum_{\dot{\sigma} \in G_x} \abs{\xi(\sigma)}^2 < \infty}.$$
We obtain a representation $\pi_x$ of $C_c(G,\Sigma)$ on $\scH_x$ by setting
$$
  \pi_x(f)(\xi)(\sigma) = \sum_{\dot{\tau} \in G_x} f(\sigma \tau^{-1}) \xi(\tau).
$$
By definition, $C^*_r(G,\Sigma)$ is the completion of $C_c(G,\Sigma)$ under the norm $\norm{f} = \sup_{x \in X} \norm{\pi_x(f)}$.

\subsection{From Cartan pairs to twisted groupoids}

Let us recall Renault's definition and characterization of Cartan subalgebras.
\setlength{\parindent}{0.5cm} \setlength{\parskip}{0cm}

\bdefin[{\cite[Definition~5.1]{R}}]
\label{CartanSubalgebra}
A sub-C*-algebra $B$ of a C*-algebra $A$ is called a Cartan subalgebra if
\begin{enumerate}
\item[(i)] $B$ contains an approximate identity of $A$;
\item[(ii)] $B$ is maximal abelian;
\item[(iii)] $B$ is regular, i.e., $N_A(B) \defeq \menge{n \in A}{n B n^* \subseteq B \ \text{and} \ n^* B n \subseteq B}$ generates $A$ as a C*-algebra;
\item[(iv)] there exists a faithful conditional expectation $P:A \onto B$.
\end{enumerate}
\edefin

\btheo[{\cite[Theorem~5.2 and Theorem~5.9]{R}}]
\label{Renault}
Cartan pairs $(A,B)$, where $A$ is a separable C*-algebra, are precisely of the form $(C^*_r(G,\Sigma),C_0(G^{(0)}))$, where $(G,\Sigma)$ is a twisted \etale{} Hausdorff locally compact second countable topologically principal groupoid.
\etheo

\bremark
\label{A,B-->GPD}
Let us briefly explain how $(G,\Sigma)$ is constructed out of $(A,B)$. Set $X \defeq \Spec(B)$. By \cite[Proposition~4.7]{R}, for every $n \in N_A(B)$ there exists a partial homeomorphism $\alpha_n: \: \dom(n) \to \ran(n)$, with $\dom(n) = \menge{x \in X}{n^*n(x) > 0}$ and $\ran(n) = \menge{x \in X}{nn^*(x) > 0}$, such that $n^*bn(x) = b(\alpha_n(x))n^*n(x)$ for all $b \in B$ and $x \in \dom(n)$. Here we use the canonical identification $B \cong C_0(X)$. Define the pseudogroup $\cG(B) \defeq \menge{\alpha_n}{n \in N_A(B)}$ on $X$. Let $G(B)$ be the groupoid of germs of $\cG(B)$, $G(B) = \menge{[x,\alpha_n,y]}{n \in N_A(B), \, y \in \dom(n), \, x = \alpha_n(y)}$. Here $[x,\alpha_n,y] = [x,\alpha_{n'},y]$ if there exists an open neighbourhood $V$ of $y$ in $X$ such that $\alpha_n \vert_V = \alpha_{n'} \vert_V$. To describe the twist, set $D \defeq \menge{(x,n,y) \in X \times N_A(B) \times X}{y \in \dom(n), \, x = \alpha_n(y)}$ and define $\Sigma(B) \defeq D / {}_{\sim}$, where $(x,n,y) \sim (x',n',y')$ if $y=y'$ and there exist $b,b' \in B$ with $b(y),b'(y) > 0$ and $nb=n'b'$. $G(B)$ and $\Sigma(B)$ are topological groupoids (see \cite{R} for details), and the canonical homomorphism $\Sigma(B) \to G(B), \, [x,n,y] \ma [x,\alpha_n,y]$, yields the central extension
$$
\Tz \times X \rightarrowtail \Sigma(B) \twoheadrightarrow G(B).
$$
In the following, we will follow \cite{R} and call $G(B)$ the Weyl groupoid and $\Sigma(B)$ the Weyl twist of $(A,B)$.
\eremark

\subsection{The Rokhlin property and approximate representability}

\bdefin \cite[Definitions~1.1]{Kir2}
For a C*-algebra $A$, its sequence algebra $A_\infty$ is defined as
$$
A_\infty := \ell^\infty(\Nz, A) / c_0(\Nz, A).
$$
There is a canonical embedding of $A$ into $A_\infty$ by (representatives) of constant sequences. If $B \subseteq A_\infty$ is a sub-C*-algebra, then we also consider the relative commutant of $B$ inside $A_\infty$
$$
A_\infty \cap B' := \menge{x \in A_\infty}{xb = bx \text{ for all } b \in B}
$$
and the (corrected) central sequence algebra
$$
F_\infty(B,A) := A_\infty \cap B' / \menge{x \in A_\infty}{xb = bx = 0 \text{ for all } b  \in B}.
$$
We also denote $F_\infty(A) = F_\infty(A,A)$.
\edefin

Observe that $F_\infty(B,A) = A_\infty \cap B'$ if $A$ is unital and $B$ contains the unit of $A$. If $B$ is $\sigma$-unital, then $F_\infty(B,A)$ is unital, see \cite[Proposition~1.9(2)+(3)]{Kir2}.

\bdefin
Let $\Gamma$ be a finite group. We denote by $\sigma:\Gamma \curvearrowright C(\Gamma)$ the canonical $\Gamma$-shift, that is, $\sigma_\gamma(f)(\eta)=f(\gamma^{-1}\cdot\eta)$ for all $f \in C(\Gamma)$ and $\gamma, \eta \in \Gamma$.
\edefin

\bdefin [cf.\ {\cite[Definition~3.1]{Iz}} and {\cite[Definition~1.3]{BS}}]
Let $\Gamma$ be a finite group and $A$ a separable C*-algebra. An action $\alpha:\Gamma\curvearrowright A$ is said to have the Rokhlin property if there exists a unital and equivariant *-homomorphism
\[
(C(\Gamma),\sigma)\to (F_\infty(A),\alpha_\infty),
\]
where $\alpha_\infty:\Gamma\curvearrowright F_\infty(A)$ is the action induced by entry-wise application of $\alpha$.
\edefin

The Rokhlin property can also be spelled out concretely in terms of elements of $A$, see for example \cite[Remark~1.4]{BS}.

Next we recall the definition of approximate representability for finite abelian group actions. For a single automorphism, this notion is already present in \cite{Kis98}.

\bdefin [{\cite[Definition~3.6]{Iz}} and {\cite[Definition~4.1]{Na}}]
Let $\Gamma$ be a finite abelian group and $A$ a separable C*-algebra. An action $\alpha: \Gamma \curvearrowright A$ is said to be approximately representable if there exists a unitary group representation $u: \Gamma \to F_\infty(A^\alpha) \subseteq F_\infty(A^\alpha,A)$ such that $\alpha_\gamma(a) = u_\gamma a u_\gamma^*$ for all $\gamma \in \Gamma$ and $a \in A$.
\edefin

For finite abelian groups, actions with the Rokhlin property and approximately representable actions are dual to each other.

\btheo [{\cite[Lemma~3.8]{Iz}} and {\cite[Proposition~4.4]{Na}}] \label{thm:duality Iz Naw}
Let $\Gamma$ be a finite abelian group and $A$ a separable C*-algebra. Let $\alpha:\Gamma \curvearrowright A$ be an action. Then
\begin{itemize}
\item[1)] $\alpha$ has the Rokhlin property if and only if $\hat{\alpha}$ is approximately representable;
\item[2)] $\alpha$ is approximately representable if and only if $\hat{\alpha}$ has the Rokhlin property.
\end{itemize}
\etheo

\bremark \label{rem: unitaries approx repr}
Let $\Gamma$ be a finite abelian group and denote by $e_\gamma \in C(\Gamma)$ the minimal projection associated with $\gamma \in \Gamma$. Let $\alpha: \Gamma \curvearrowright A$ be an action with the Rokhlin property on a unital and separable C*-algebra. Let us remark for later purposes that for a given unital and equivariant *-homomorphism $\iota: (C(\Gamma),\sigma)\to (A_\infty \cap A',\alpha_\infty)$, the map 
\[
u:\hat{\Gamma} \to A_\infty \cap A', \ u_\chi = \sum_{\gamma \in \Gamma} \chi(\gamma) \iota(e_\gamma)
\]
defines a unitary representation of the dual group $\hat{\Gamma}$. As explained in \cite[Lemma~3.8]{Iz}, $\hat{\alpha}_\chi(x) = u_\chi^* x u_\chi$ for all $x \in A \rtimes_\alpha \Gamma$ and $\chi \in \hat{\Gamma}$. In particular, unitaries in $A_\infty = ((A \rtimes_\alpha \Gamma)^{\hat{\alpha}})_\infty$ witnessing approximate representability of $\hat{\alpha}$ can be chosen in $\iota(C(\Gamma))$.
\eremark

\section{Inductive limits and Cartan pairs}
\label{Sec:indlim}

Let $(G,\Sigma)$ and $(H,T)$ be twisted groupoids. A twisted groupoid homomorphism is a continuous groupoid homomorphism $p: \: T \to \Sigma$ which is $\Tz$-equivariant, i.e., $p(z \tau) = z p(\tau)$ for all $z \in \Tz$ and $\tau \in T$. Let $\pi_H: \: T \onto H$ and $\pi_G: \: \Sigma \onto G$ be the canonical projections. $p$ induces a map $\dot{p}: \: H \to G, \, \dot{\tau} \ma \pi_G(p(\tau))$. In other words, $\dot{p}$ is uniquely determined by the commutative diagram
\begin{align}
\label{CD_pdotp}
  \xymatrix{
  \Sigma \ar[d]_{\pi_G} & \ar[l]_{p} T \ar[d]^{\pi_H}
  \\
  G & \ar[l]_{\dot{p}} H
  }
\end{align}

\blemma
\label{Lem:p-proper}
In the situation above, $p$ is proper if and only if $\dot{p}$ is proper.
\elemma
\bproof
This is an easy consequence of the observation that $\pi_G$ and $\pi_H$ are proper and commutativity of \eqref{CD_pdotp}.
\eproof

\blemma
\label{Lem:proper--hom}
Let $(G,\Sigma)$ and $(H,T)$ be twisted \etale{} groupoids. Write $X \defeq G^{(0)}$ and $Y \defeq H^{(0)}$. Let $p: \: T \to \Sigma$ be a surjective twisted groupoid homomorphism. Assume that $p$ is proper, and that for every $y \in Y$, $H_y \to G_{p(y)}, \, \eta \ma \dot{p}(\eta)$ is a bijection. Then $C_c(G,\Sigma) \to C_c(H,T), \, f \ma f \circ p$ extends to an isometric homomorphism $p^*: \: C^*_r(G,\Sigma) \to C^*_r(H,T)$.
\elemma
\bproof
By Lemma~\ref{Lem:p-proper} and because $p$ is proper, $\dot{p}$ is proper. Thus given $f \in C_c(G,\Sigma)$, we have that $\supp(f \circ p) = \dot{p}^{-1}(\supp(f))$ is compact because $\supp(f)$ is compact. Hence $f \circ p$ lies in $C_c(H,T)$. As $p$ is a groupoid homomorphism, it is clear that $(f \circ p)^* = f^* \circ p$. To see multiplicativity, take $f, g \in C_c(G,\Sigma)$ and compute
$$
  (f * g)(p(\tau)) = \sum_{\dot{\sigma} \in G_{s(p(\tau))}} f(p(\tau) \sigma^{-1}) g(\sigma) = \sum_{\dot{\upsilon} \in H_{s(\tau)}} f(p(\tau)p(\upsilon)^{-1}) g(p(\upsilon)) = ((f \circ p) * (g \circ p))(\tau).
$$
Here we used that $H_{s(\tau)} \to G_{p(s(\tau))} = G_{s(p(\tau))}, \, \eta \ma \dot{p}(\eta)$ is a bijection.

To see that $f \ma f \circ p$ is isometric, take $y \in Y$ and let $x \defeq p(y)$. Write $\scH_y$ and $\scH_x$ for the Hilbert spaces as in \S~\ref{ss:twistGPD,CSTAR}. Because $H_y \to G_x, \, \eta \ma \dot{p}(\eta)$ is a bijection, $U_p: \: \scH_x \to \scH_y, \, \xi \ma \xi \circ p$ is a unitary. For every $\xi \in \scH_x$, we have
\begin{align*}
  U_p(\pi_x(f)(\xi))(\tau) &= \pi_x(f)(\xi)(p(\tau)) = \sum_{\dot{\sigma} \in G_x} f(p(\tau) \sigma^{-1}) \xi(\sigma)\\
  &= \sum_{\dot{\upsilon} \in H_y} f(p(\tau)p(\upsilon)^{-1}) \xi(p(\upsilon)) = \sum_{\dot{\upsilon} \in H_y} f(p(\tau \upsilon^{-1})) \xi(p(\upsilon)) = \pi_y(f \circ p)(U_p(\xi))(\tau).
\end{align*}
Hence $\pi_y(f \circ p) \sim_u \pi_x(f)$, and thus,
$$
  \bigoplus_{y \in Y} \pi_y(f \circ p) \sim_u \bigoplus_{x \in X} \bigoplus_{y \in \dot{p}^{-1}(x)} \pi_x(f) \sim_u \bigoplus_{x \in X} \rukl{I_{\ell^2(\dot{p}^{-1}(x))} \otimes \pi_x(f)},
$$
where $I$ denotes the identity operator. This shows that
$$
  \norm{f \circ p}_{C^*_r(H,T)} = \norm{\bigoplus_{y \in Y} \pi_y(f \circ p)} = \norm{\bigoplus_{x \in X} \rukl{I_{\ell^2(\dot{p}^{-1}(x))} \otimes \pi_x(f)}} = \norm{\bigoplus_{x \in X} \pi_x(f)} = \norm{f}_{C^*_r(G,\Sigma)}.
$$
\eproof

Now let $(H,T)$ and $(G,\Sigma)$ be two twisted \etale{} groupoids. Let $\imath: \: T \to \Sigma$ be a twisted groupoid homomorphism, and $i: \: H \to G$ the induced groupoid homomorphism.
\blemma
\label{Lem:i-openimage}
In the situation above, $\imath(T)$ is an open subgroupoid of $\Sigma$ if and only if $i(H)$ is an open subgroupoid of $G$.
\elemma
\bproof
As $\pi_G$ is open and continuous, this follows from $\imath(T) = \pi_G^{-1}(i(H))$, which is easy to see.
\eproof
The next result is the analogue of \cite[Proposition~1.9]{Phi05} for twisted groupoids.
\blemma
\label{Lem:i:opensubgpd-hom}
Let $(H,T)$ and $(G,\Sigma)$ be twisted \etale{} groupoids. Let $\imath: \: T \into \Sigma$ be an injective twisted groupoid homomorphism such that $\imath(T)$ is an open subgroupoid in $\Sigma$. Then $C_c(H,T) \to C_c(G,\Sigma), \, f \ma 1_T \cdot f$ extends to an isometric homomorphism $\imath_*: \: C^*_r(H,T) \to C^*_r(G,\Sigma)$.
\elemma
\setlength{\parindent}{0cm} \setlength{\parskip}{0cm}

Here we identify $T$ with an open subgroupoid of $\Sigma$ via $\iota$. $1_T$ is the characteristic function of $T$, and $1_T \cdot f$ is the extension of $f$ to the function on $\Sigma$ such that $(1_T \cdot f)(\sigma) = 0$ if $\sigma \notin T$ and $(1_T \cdot f)(\sigma) = f(\sigma)$ if $\sigma \in T$.
\setlength{\parindent}{0.5cm} \setlength{\parskip}{0cm}

\bproof
Let $i: \: H \to G$ be the groupoid homomorphism induced by $\imath$. As $\imath$ is injective, $i$ must be so as well. Let us identify $H$ and $T$ with open subgroupoids of $G$ and $\Sigma$ via $i$ and $\imath$, respectively. The image of $i$ is open because of Lemma~\ref{Lem:i-openimage}. As $T$ is open in $\Sigma$, $f \ma 1_T \cdot f$ is well-defined (the property of having compact support is clearly preserved). It is moreover easy to see that $f \ma 1_T \cdot f$ is a *-homomorphism. It remains to prove that $f \ma 1_T \cdot f$ is isometric.

Fix $x \in X \defeq G^{(0)}$. For $\alpha_1, \alpha_2 \in \Sigma_x$, define $\alpha_1 \sim \alpha_2$ if $\alpha_1 \alpha_2^{-1} \in T$. Similarly, for $\gamma_1, \gamma_2 \in G_x$, define $\gamma_1 \sim \gamma_2$ if $\gamma_1 \gamma_2^{-1} \in H$. It is easy to see that, for all $\alpha_1, \alpha_2 \in \Sigma_x$, $\alpha_1 \sim \alpha_2$ if and only if $\dot{\alpha}_1 \sim \dot{\alpha}_2$. Moreover, for $\alpha \in \Sigma_x$, we have $\alpha \sim \alpha$ if and only if $r(\alpha) = \alpha \alpha^{-1} \in T$ if and only if $r(\alpha) \in Y \defeq H^{(0)}$. ($r$ here is the range map of $\Sigma$.) Define $\Rho_x \defeq \menge{\alpha \in \Sigma_x}{\alpha \sim \alpha} = r^{-1}(Y) \cap \Sigma_x$. Restricted to $\Rho_x$, $\sim$ is an equivalence relation. For $\alpha \in \Rho_x$, consider its equivalence class $[\alpha] \defeq \menge{\beta \in \Sigma_x}{\beta \sim \alpha}$. We have a $\Tz$-equivariant homeomorphism $[\alpha] \to T_{r(\alpha)}, \, \beta \ma \beta \alpha^{-1}$, whose inverse is given by $\tau \alpha \mafr \tau$. Let $\scH_x$ be the Hilbert space as in \S~\ref{ss:twistGPD,CSTAR}. Given $\alpha \in \Rho_x$, define $\scH_{[\alpha]} \defeq \menge{\xi \in \scH_x}{\xi(\beta) = 0 \ {\rm for} \ {\rm all} \ \beta \notin [\alpha]}$. Also, let $\scH_{\Rho_x^{\mathsf{c}}} \defeq \menge{\xi \in \scH_x}{\xi(\sigma) = 0 \ {\rm for} \ {\rm all} \ \sigma \in \Rho_x}$. Then $\scH_x = \rukl{\bigoplus_{[\alpha] \in {\rm P}_x / {}_{\sim}} \scH_{[\alpha]}} \oplus \scH_{\Rho_x^{\mathsf{c}}}$. For $y \in Y$, let $\scH_y$ be the Hilbert space as in \S~\ref{ss:twistGPD,CSTAR}. For $\alpha \in \Rho_x$, consider the unitary $U_{\alpha}: \: \scH_{[\alpha]} \to \scH_{r(\alpha)}, \xi \ma \xi(\sqcup \alpha)$. Write $[\dot{\alpha}] \defeq \menge{\gamma \in G_x}{\gamma \sim \dot{\alpha}}$. Given $\xi \in \scH_{[\alpha]}$ and $\beta \in \Sigma_x$,
$$
  \pi_x(1_T \cdot f)(\xi)(\beta) = \sum_{\dot{\sigma} \in G_x} (1_T \cdot f)(\beta \sigma^{-1}) \xi(\sigma) = \sum_{\dot{\sigma} \in [\dot{\alpha}]} (1_T \cdot f)(\beta \sigma^{-1}) \xi(\sigma)
$$
vanishes unless $\beta \sigma^{-1} \in T$, i.e., $\beta \sim \sigma \sim \alpha$. Hence $\pi_x(1_T \cdot f)(\scH_{[\alpha]}) \subseteq \scH_{[\alpha]}$. We have
\begin{align*}
  \pi_{r(\alpha)}(f)(U_{\alpha} \xi)(\tau) &= \sum_{\dot{\upsilon} \in H_{r(\alpha)}} f(\tau \upsilon^{-1}) \xi(\upsilon \alpha) = \sum_{\dot{\upsilon} \in H_{r(\alpha)}} f((\tau \alpha) (\upsilon \alpha)^{-1}) \xi(\upsilon \alpha)\\
  &= \sum_{\dot{\beta} \in [\dot{\alpha}]} f(\tau \alpha \beta^{-1}) \xi(\beta) = \sum_{\dot{\beta} \in G_x} (1_T \cdot f)(\tau \alpha \beta^{-1}) \xi(\beta) = \pi_x(1_T \cdot f)(\xi)(\tau \alpha)\\
  &= U_{\alpha}(\pi_x(1_T \cdot f)(\xi))(\tau).
\end{align*}
This shows that $\pi_x(1_T \cdot f) \vert_{\scH_{[\alpha]}} \sim_u \pi_{r(\alpha)}(f)$. It is clear that $\pi_x(1_T \cdot f) \vert_{\scH_{\Rho_x^{\mathsf{c}}}} = 0_{\Rho_x^{\mathsf{c}}}$. Thus $\pi_x(1_T \cdot f) \sim_u \rukl{\bigoplus_{[\alpha] \in \Rho_x / {}_{\sim}} \pi_{r(\alpha)}(f)} \oplus 0_{\Rho_x^{\mathsf{c}}}$. Hence
\begin{align*}
  \norm{1_T \cdot f}_{C^*_r(G,\Sigma)} &= \norm{\bigoplus_{x \in X} \pi_x(1_T \cdot f)} = \norm{\bigoplus_{x \in X} \rukl{\rukl{\bigoplus_{[\alpha] \in \Rho_x / {}_{\sim}} \pi_{r(\alpha)}(f)} \oplus 0_{\Rho_x^{\mathsf{c}}}}}\\
  &= \sup \menge{\norm{\pi_{r(\alpha)}(f)}}{[\alpha] \in \Rho_x / {}_{\sim}, \, x \in X} = \sup \menge{\norm{\pi_y(f)}}{y \in Y} = \norm{f}_{C^*_r(H,T)}.
\end{align*}
\eproof

Recall that for a locally compact Hausdorff space $X$, a sequence $(f_n)_n \in C_c(X)$ is said to converge to $f \in C_c(X)$ in the inductive limit topology if there exists a compact subset $K \subseteq X$ such that $\supp(f_n), \supp(f) \subseteq K$ for all $n$ and $f_n$ converges to $f$ on $K$ uniformly. 

\blemma
\label{Lem:plim-dense}
Let $G_j$ be locally compact Hausdorff second countable topological groupoids and $p_j: \: G_{j+1} \to G_j$ continuous proper groupoid homomorphisms. Form the inverse limit $\bG \defeq \plim_j \gekl{G_j; p_j}$ and let $p_j^{\infty}:\bG \to G_j$ be the canonical projections.

Then $\bG$ becomes a locally compact Hausdorff second countable topological groupoid such that all $p_j^{\infty}$ are proper and $\bigcup_j (p_j^{\infty})^* (C_c(G_j))$ is dense in $C_c(\bG)$ in the inductive limit topology.
\elemma
\bproof
By definition,
$$
  \bG = \plim_j \gekl{G_j; p_j}
  = \menge{(\gamma_j) \in \prod_j G_j}{p_j(\gamma_{j+1}) = \gamma_j} \subseteq \prod_j G_j,
$$
and the topology on $\bG$ is given by the subspace topology from $\prod_j G_j$. It is clear that $\bG$ is Hausdorff and second countable. To see that $\bG$ is locally compact, take an arbitrary element ${\bar \gamma} = (\gamma_j)$ in $\bG$. Choose a compact neighbourhood $U_0$ of $\gamma_0 \in G_0$. Define recursively $U_{j+1} \defeq p_j^{-1}(U_j)$. Since all $p_j$ are proper, all the $U_j$ are compact. Hence
$$
  (p_0^\infty)^{-1}(U_0) = \menge{(\eta_j) \in \bG}{\eta_j \in U_j} \subseteq \prod_j U_j
$$
is a compact neighbourhood of $\bar{\gamma}$. Therefore, with component-wise groupoid operations $\bG$ becomes a locally compact Hausdorff second countable topological groupoid.

Now let us show that $p^\infty_j$ is proper. Take $K_j \subseteq G_j$ compact. Define recursively $K_{i+1} \defeq p_i^{-1}(K_i)$ if $i \geq j$ and $K_{i-1} \defeq p_{i-1}(K_i)$ if $i \leq j$. This gives us compact subsets $K_i \subseteq G_i$, for all $i$, and we obviously have that $(p^\infty_j)^{-1}(K_j) = \menge{(\gamma_i) \in \bG}{\gamma_i \in K_i \ {\rm for} \ {\rm all} \ i}$ is compact.

The proof that $\bigcup_j (p_j^{\infty})^* (C_c(G_j))$ is dense in $C_c(\bG)$ in the inductive limit topology is a straightforward partition of unity argument, which we leave to the reader.
\eproof

For $k = 1, 2, \dotsc$, let $(A_k,B_k) = (C^*_r(G_k,\Sigma_k),C_0(X_k))$ for twisted \etale{} groupoids $(G_k,\Sigma_k)$ given by 
$$\Tz \times X_k \overset{\iota_k}{\rightarrowtail} \Sigma_k \overset{\pi_k}{\twoheadrightarrow} G_k,$$
where $X_k = G_k^{(0)}$. Assume that, for all $k$, $\vp_k: \: A_k \to A_{k+1}$ is an injective *-homomorphism with $\vp_k(B_k) \subseteq B_{k+1}$. Further assume that, for all $k$, $(H_k,T_k)$ is a twisted \etale{} groupoid such that $H_k$ is an open subgroupoid of $G_{k+1}$ and $T_k = \pi_{k+1}^{-1}(H_k)$. For all $k$, let $\imath_k: \: T_k \to \Sigma_{k+1}$ and $i_k: \: H_k \to G_{k+1}$ be the canonical inclusions. For all $k$, let $p_k: \: T_k \to \Sigma_k$ be a twisted groupoid homomorphism. Assume that, for all $k$, $p_k$ is surjective, proper, and that for all $y \in Y_k \defeq H_k^{(0)}$, $(H_k)_y \to (G_k)_{p_k(y)}, \, \eta \ma \dot{p}_k(\eta)$ is a bijection. Further assume that for all $k$, $\vp_k = (\imath_k)_* \circ p_k^*$. Suppose that the following two conditions are satisfied:
\begin{enumerate}
\item[(LT)] For every open subset $U \subseteq G_k$ and continuous section $\rho: \: U \to \Sigma_k$ for $\pi_k$, there is a continuous section $\ti{\rho}: \: \dot{p}_k^{-1}(U) \to T_k$ for $\pi_{k+1} \vert_{T_k}$, such that
\begin{align*}
  \xymatrix{
  \Sigma_k & \ar[l]_{p_k} T_k
  \\
  U \ar[u]^{\rho} & \ar[l]_{\dot{p}_k} \dot{p}_k^{-1}(U) \ar[u]_{\ti{\rho}}
  }
\end{align*}
commutes.
\item[(E)] For every open subset $U \subseteq X_k$ and continuous section $t: \: U \to G_k$ for the source map of $G_k$, where $t$ has open image, there is a continuous section $\ti{t}: \: \dot{p}_k^{-1}(U) \to H_k$ for the source map of $H_k$, where $\ti{t}$ has open image, such that
\begin{align*}
  \xymatrix{
  G_k & \ar[l]_{\dot{p}_k} H_k
  \\
  U \ar[u]^{t} & \ar[l]_{\dot{p}_k} \dot{p}_k^{-1}(U) \ar[u]_{\ti{t}}
  }
\end{align*}
commutes.
\end{enumerate}
For all $k$, define recursively for all $j = 0, 1, \dotsc$: $\Sigma_{k,0} \defeq \Sigma_k$ and $\Sigma_{k,j+1} \defeq p_{k+j}^{-1}(\Sigma_{k,j}) \subseteq T_{k+j}$. Similarly, for all $k$, define recursively for all $j = 0, 1, \dotsc$: $G_{k,0} \defeq G_k$ and $G_{k,j+1} \defeq \dot{p}_{k+j}^{-1}(G_{k,j}) \subseteq H_{k+j}$. Set $\bSigma_k \defeq \plim_j \gekl{\Sigma_{k,j}; p_{k+j}}$ and $\bG_k \defeq \plim_j \gekl{G_{k,j}; \dot{p}_{k+j}}$. In this situation, we have
\btheo
\label{indlim--GPD}
{$ $}
\begin{enumerate}
\item[(a)] For all $k$, $(\bG_k, \bSigma_k)$ is a twisted \etale{} groupoid. Moreover, the inclusion maps $\imath_k$ and $i_k$ induce, for all $k$, an injective twisted groupoid homomorphism ${\bar \imath}_k: \: \bSigma_k \to \bSigma_{k+1}$ and an injective continuous groupoid homomorphism ${\bar i}_k: \: \bG_k \to \bG_{k+1}$ with open image. Setting $\bSigma \defeq \ilim_k \gekl{\bSigma_k; {\bar \imath}_k}$ and $\bG \defeq \ilim_k \gekl{\bG_k; {\bar i}_k}$, $(\bG, \bSigma)$ is a twisted \etale{} groupoid.
\item[(b)] We have $\ilim_k \gekl{A_k; \vp_k} \cong C^*_r(\bG, \bSigma)$ and $\ilim_k \gekl{B_k; \vp_k} \cong C_0(\bG^{(0)})$.
\item[(c)] If for all $k$, $G_k$ is topologically principal, then also $\bG_k$ is topologically principal for all $k$, and $\bG$ is topologically principal. In that case, $\rukl{\ilim_k \gekl{A_k; \vp_k}, \ilim_k \gekl{B_k; \vp_k}} \cong \rukl{C^*_r(\bG, \bSigma), C_0(\bG^{(0)})}$ is a Cartan pair.
\end{enumerate}
\etheo
\bproof
Let us first prove (a) and (b). By definition,
$$
  \bSigma_k = \plim_j \gekl{\Sigma_{k,j}; p_{k+j}}
  = \menge{(\sigma_j) \in \prod_j \Sigma_{k,j}}{p_{k+j}(\sigma_{j+1}) = \sigma_j} \subseteq \prod_j \Sigma_{k,j},
$$
and the topology on $\bSigma_k$ is given by the subspace topology from $\prod_j \Sigma_{k,j}$. Let $p^k_j: \: \bSigma_k \onto \Sigma_{k,j}$ be the canonical projection. Similarly,
$$
  \bG_k = \plim_j \gekl{G_{k,j}; \dot{p}_{k+j}}
  = \menge{(\gamma_j) \in \prod_j G_{k,j}}{\dot{p}_{k+j}(\gamma_{j+1}) = \gamma_j} \subseteq \prod_j G_{k,j},
$$
with the subspace topology from $\prod_j G$. Let $\dot{p}^k_j: \: \bG_k \onto G_{k,j}$ be the canonical projection.

It follows from Lemma~\ref{Lem:p-proper} and Lemma~\ref{Lem:plim-dense} that, with component-wise groupoid operations, both $\bG_k$ and $\bSigma_k$ become locally compact Hausdorff second countable topological groupoids.

Let $X_{k,j} \defeq G_{k,j}^{(0)}$ and $\bX_k \defeq \bG_k^{(0)}$. Clearly, we have
$$
  \bX_k = \menge{(\gamma_j) \in \bG_k}{\gamma_j \in X_{k,j} \ {\rm for} \ {\rm all} \ j}.
$$
Hence
\begin{align}
\label{barG_k^0}
  \bX_k = \plim_j \gekl{X_{k,j}; \dot{p}_{k+j}}.
\end{align}
To show that $\bG_k$ is \etale, we have to prove that the source map ${\bar s}: \: \bG_k \to \bX_k, \, (\gamma_j) \ma (s(\gamma_j))$ is a local homeomorphism. Equivalently, we show that for any given ${\bar \gamma} = (\gamma_j) \in \bG_k$, there exists an open subset $\bar{U} \subseteq \bX_k$ containing ${\bar x} = {\bar s}({\bar \gamma}) = (x_j)$ and a continuous section $\bar{t}: \bar{U} \to \bG_k$ with open image for $\bar{s}$ such that $\bar{t}(\bar{x}) = {\bar \gamma}$. Let $U \subseteq X_k$ be an open neighbourhood of $x_0$ and $t_0: \: U \to G_k$ a continuous section with open image for $s$ (the source map of $G_k$), with $t_0(x_0) = \gamma_0$. Define recursively $U_{j+1} \defeq \dot{p}_{k+j}^{-1}(U_j)$. By condition (E), we can recursively find continuous sections $t_{j+1}: \: U_{j+1} \to G_{k,j+1}$ with open image for the source map of $G_{k,j+1}$ such that
\begin{align*}
  \xymatrix{
  G_{k,j} & \ar[l]_{\dot{p}_{k+j}} G_{k,j+1}
  \\
  U_j \ar[u]^{t_j} & \ar[l]_{\dot{p}_{k+j}} U_{j+1} \ar[u]_{t_{j+1}}
  }
\end{align*}
commutes. Define $\bar U : = (\dot{p}_0^k)^{-1}(U_0)$. Then ${\bar t}: \: \bar U  \to \bG_k, \, (x_j) \ma (t_j(x_j))$ is a continuous section with open image for ${\bar s}$. In addition, we must have ${\bar t}({\bar x}) = {\bar \gamma}$, i.e., $(t_j(x_j)) = (\gamma_j)$. This follows inductively from $t_0(x_0) = \gamma_0$, $\dot{p}_{k+j}(t_{j+1}(x_{j+1})) = t_j(x_j)$, $\dot{p}_{k+j}(\gamma_{j+1}) = \gamma_j$ and the assumption that $\dot{p}_{k+j} \vert_{(G_{k,j+1})_{x_{j+1}}}:(G_{k,j+1})_{x_{j+1}} \to (G_{k,j})_{x_j}$ is bijective. Hence $\bG_k$ is \etale, for all $k$.

Define $\Tz \curvearrowright \bSigma_k$ by $z (\sigma_j) = (z \sigma_j)$ and ${\bar \pi}_k: \: \bSigma_k \to \bG_k, \, (\sigma_j) \ma (\pi_{k+j}(\sigma_j))$. Obviously, the $\Tz$-action is continuous and free. Assume that for $(\sigma_j) \in \bSigma_k$, ${\bar \pi}_k(\sigma_j) = (x_j) \in \bX_k$. Then we know that there are $z_j \in \Tz$ with $\sigma_j = z_j x_j$ for all $j$. Let $z \defeq z_0$. Then it follows from
$$
  z_{j+1} x_j = z_{j+1} \dot{p}_{k+j}(x_{j+1}) = p_{k+j}(\sigma_{j+1}) = \sigma_j = z_j x_j
$$
that $z_{j+1} = z_j$, and hence $z_j = z$ for all $j$. Thus the canonical projection $\bSigma_k / \Tz \to \bG_k$ induced by ${\bar \pi}_k$ is continuous and injective, which shows that $\bSigma_k / \Tz$ is Hausdorff. Let us show that ${\bar \pi}_k: \: \bSigma_k \to \bG_k$ is surjective and defines a locally trivial $\Tz$-bundle. Let ${\bar \gamma} = (\gamma_j) \in \bG_k$. Let $U_0 \subseteq G_k$ be an open neighbourhood of $\gamma_0$ and $\rho_0: \: U_0 \to \Sigma_k$ a continuous section for $\pi_k$. $U_0$ and $\rho_0$ exist because $\pi_k: \: \Sigma_k \onto G_k$ is a locally trivial $\Tz$-bundle. Define recursively $U_{j+1} \defeq \dot{p}_{k+j}^{-1}(U_j)$. By condition (LT), we can recursively find continuous sections $\rho_{j+1}: \: U_{j+1} \to \Sigma_{k,j+1}$ for $\pi_{k+j+1}$ such that
\begin{align*}
  \xymatrix{
  \Sigma_{k,j} & \ar[l]_{p_{k+j}} \Sigma_{k,j+1}
  \\
  U_j \ar[u]^{\rho_j} & \ar[l]_{\dot{p}_{k+j}} U_{j+1} \ar[u]_{\rho_{j+1}}
  }
\end{align*}
commutes. Define $\bar U := (\dot{p}^k_0)^{-1}(U_0)$. Then ${\bar \rho}: \: \bar U   \to \bSigma_k, \, (\gamma_j) \ma (\rho_j(\gamma_j))$ is a continuous section for ${\bar \pi}_k$. In particular, ${\bar \pi}_k({\bar \rho}({\bar \gamma})) = {\bar \gamma}$, and we see that ${\bar \pi}_k$ is surjective. We now obtain homeomorphisms 
\begin{align}
\label{cd:homeo loc triv pi_{k+j}}
\Tz \times U_j \overset{\cong}{\lori} \pi_{k+j}^{-1}(U_j), \, (z,\gamma_j) \ma z \rho_j(\gamma_j) \text{ for all } z \in \Tz,\ \gamma_j \in U_j
\end{align}
with inverse $(z_j(\sigma),\dot{\sigma}) \mafr \sigma$. It follows that 
\begin{align}
\label{cd:homeo loc triv}
\Tz \times {\bar U} \overset{\cong}{\lori} {\bar \pi}_k^{-1}({\bar U}), \, (z,\bar \gamma) \ma z {\bar \rho}(\bar \gamma) \text{ for all } z \in \Tz,\ \bar \gamma \in \bar U
\end{align}
defines a homeomorphism with inverse $(z(\sigma),\dot{\sigma}) \mafr \sigma$. Hence, ${\bar \pi}_k$ indeed defines a locally trivial principal $\Tz$-bundle. Thus, for all $k$, $(\bG_k,\bSigma_k)$ is a twisted \etale{} groupoid.

By Lemma~\ref{Lem:plim-dense} and Lemma~\ref{Lem:p-proper}, all $p^k_j$ and $\dot{p}^k_j$ are proper. 
Moreover, $p^k_j$ is surjective since $p_{k+j}$ is surjective for all $j$. Now let us show that, for every ${\bar x} = (x_j) \in \bX_k$, the restriction of $\dot{p}^k_j$ to $(\bG_k)_{\bar x}$ is bijective. It suffices to show that $(\gamma_i) \in (\bG_k)_{\bar x}$ is uniquely determined by $\gamma_j = \dot{p}^k_j((\gamma_i))$. This now follows inductively as, for $i \geq j$, $\gamma_{i+1}$ is determined by $\gamma_i$ because of $\dot{p}_{k+i}(\gamma_{i+1}) = \gamma_i$ and $\dot{p}_{k+i}$ induces a bijection $(G_{k,i+1})_{x_{i+1}} \to (G_{k,i})_{x_i}$, while for $i \leq j$, $\dot{p}_{k+i-1}(\gamma_i) = \gamma_{i-1}$ shows that $\gamma_i$ determines $\gamma_{i-1}$.

Thus Lemma~\ref{Lem:proper--hom} tells us that $C_c(G_{k,j},\Sigma_{k,j}) \to C_c(\bG_k,\bSigma_k), \, f \ma f \circ p_j^k$ extends to an isometric homomorphism $(p_j^k)^*: \: C^*_r(G_{k,j},\Sigma_{k,j}) \to C^*_r(\bG_k,\bSigma_k)$. Let us now show that this collection $(p_j^k)^*$ gives rise to an isomorphism $\ilim_j \gekl{C^*_r(G_{k,j},\Sigma_{k,j}); p_{k+j}^*} \overset{\cong}{\lori} C^*_r(\bG_k,\bSigma_k)$ which identifies $\ilim_j \gekl{C_0(X_{k,j}); p_{k+j}^*}$ with $C_0(\bX_k)$.

We have $(p_{j+1}^k)^* \circ (p_{k+j})^* = (p_{k+j} \circ p_{j+1}^k)^* = (p_j^k)^*$, so that the $(p_j^k)^*$ indeed induce a homomorphism $\ilim_j (p_j^k)^*: \: \ilim_j \gekl{C^*_r(G_{k,j},\Sigma_{k,j}); p_{k+j}^*} \to C^*_r(\bG_k,\bSigma_k)$. As all the $(p_j^k)^*$ are isometric, $\ilim_j (p_j^k)^*$ must be isometric as well. Hence it suffices to show that $\ilim_j (p_j^k)^*$ has dense image. For that purpose, it is enough to show that for every $f \in C_c(\bG_k,\bSigma_k)$, there exists a sequence $f_m \in \bigcup_j (p_j^k)^*(C_c(G_{k,j},\Sigma_{k,j}))$ such that $f_m \lori_m f$ in the inductive limit topology. The latter means that there is a compact subset $K \subseteq \bG_k$ such that $\supp(f_m) \subseteq K$ for all $m$, and we have $\lim_{m \to \infty} \norm{f - f_m}_{\infty} = 0$.

Using a partition of unity, we may assume that $f \in C_c(\bG_k,\bSigma_k)$ has $\supp(f)$ contained in an open subset ${\bar U}$ of the form ${\bar U} = (\dot{p}_0^k)^{-1}(U_0)$, where $U_0 \subseteq G_k$ is an open subset for which there exists a continuous section $\rho_0: \: U_0 \to \Sigma_k$ for $\pi_k$. As above, using condition (LT), find a continuous section ${\bar \rho}: \: {\bar U} = (\dot{p}_0^k)^{-1}(U_0) \to \bSigma_k$ for ${\bar \pi}_k$. We have ${\bar \rho} = (\rho_j)$, where $\rho_j: \: U_j \to \Sigma_{k,j}$ are continuous sections for $\pi_{k+j}$, such that
\begin{align}
\label{cd:rhorho_j}
  \xymatrix{
  \Sigma_{k,j} & \ar[l]_{p^k_j} \bSigma_k
  \\
  U_j \ar[u]^{\rho_j} & \ar[l]_{\dot{p}^k_j} {\bar U} \ar[u]_{{\bar \rho}}
  }
\end{align}
commutes. Recall the homeomorphisms $\pi_{k+j}^{-1}(U_j) \overset{\cong}{\lori} \Tz \times U_j$, $\sigma \mapsto (z_j(\sigma),\dot{\sigma})$ and ${\bar \pi}_k^{-1}({\bar U}) \overset{\cong}{\lori} \Tz \times {\bar U}$, $\sigma \mapsto (z(\sigma),\dot{\sigma})$ given as in \eqref{cd:homeo loc triv pi_{k+j}} and \eqref{cd:homeo loc triv}, respectively. Because \eqref{cd:rhorho_j} commutes, we must have 
\begin{align}
\label{zsigma}
z(\sigma) = z_j(p^k_j(\sigma)) \ {\rm for} \ {\rm all} \ \sigma \in {\bar \pi}_k^{-1}({\bar U}).
\end{align}
Now, for our given $f \in C_c(\bG_k,\bSigma_k)$ with $\supp(f) \subseteq {\bar U}$, we find $\dot{f} \in C_c(\bG_k)$ with $\supp(\dot{f}) \subseteq {\bar U}$ such that $f(\sigma) = \overline{z(\sigma)} \dot{f}(\dot{\sigma})$ for all $\sigma \in {\bar \pi}_k^{-1}({\bar U}) \cong \Tz \times {\bar U}$. By Lemma~\ref{Lem:plim-dense}, we can find a fixed compact subset $K \subseteq \bG_k$ with $K \subseteq {\bar U}$ and for all $\ve > 0$ a function $\dot{g} \in C_c(G_{k,j})$ satisfying $\supp(\dot{g} \circ \dot{p}^k_j) \subseteq K$ and $\norm{\dot{f} - \dot{g} \circ \dot{p}^k_j}_{\infty} < \ve$. As $\dot{p}_j^k(K) \subseteq U_j$, we can arrange $\supp(\dot{g}) \subseteq U_j$ by possibly multiplying $\dot{g}$ with some $\varphi \in C_c(G_{k,j})$ with $\supp(\varphi) \subset U_j$ and $\varphi(\dot{\sigma}) = 1$ for all $\dot{\sigma} \in K$. We may thus define $g \in C_c(G_{k,j},\Sigma_{k,j})$ with $\supp(g) \subseteq U_j$ and $g(\sigma) = \overline{z_j(\sigma)} \dot{g}(\dot{\sigma})$ for all $\sigma \in \pi_{k+j}^{-1}(U_j) \cong \Tz \times U_j$. We have $\supp(g \circ p_j^k) = \supp(\dot{g} \circ \dot{p}_j^k) \subseteq K$ and, for all $\sigma \in \bSigma_k$,
$$
  \abs{f(\sigma) - g(p_j^k(\sigma))} = \abs{f(\sigma) - \overline{z_j(p_j^k(\sigma))} \dot{g}(\dot{p}_j^k(\dot{\sigma}))} \overset{\eqref{zsigma}}{=} \abs{\overline{z(\sigma)} \rukl{\dot{f}(\dot{\sigma}) - \dot{g}(\dot{p}_j^k(\dot{\sigma}))}} < \ve.
$$
Hence $\ilim_j (p_j^k)^*$ is indeed an isomorphism $\ilim_j \gekl{C^*_r(G_{k,j},\Sigma_{k,j}); p_{k+j}^*} \overset{\cong}{\lori} C^*_r(\bG_k,\bSigma_k)$.

Since the restriction of $\dot{p}_{k+j}$ to $(G_{k,j+1})_y$ is injective for all $y \in X_{k,j+1}$, we have $p_{k+j}^{-1}(X_{k,j}) \subseteq X_{k,j+1}$. To see this, let $\eta \in \Sigma_{k,j+1}$ with $p_{k+j}(\eta) = x \in X_{k,j}$. Let $y = s(\eta)$. Then both $\dot{\eta}$ and $y$ lie in $(G_{k,j+1})_y$, and $\dot{p}_{k+j}(\dot{\eta}) = \dot{p}_{k+j}(y) = x$. Hence $\dot{\eta} = y$. This implies $\eta = z y$ for some $z \in \Tz$, and $x = p_{k+j}(\eta) = p_{k+j}(zy) = zx$ implies $z = 1$. Thus $\eta = y \in X_{k,j+1}$. This shows $p_{k+j}^{-1}(X_{k,j}) \subseteq X_{k,j+1}$. We conclude that $\ilim_j (p_j^k)^*$ sends $\ilim_j \gekl{C_0(X_{k,j}); p_{k+j}^*}$ into $C_0(\bX_k)$. Another application of Lemma~\ref{Lem:plim-dense} yields that $\rukl{\ilim_j (p_j^k)^*} \rukl{\ilim_j \gekl{C_0(X_{k,j}); p_{k+j}^*}}$ is dense in $C_0(\bX_k)$, so that \begin{align}
\label{ilimX_kj}
  {\ilim_j (p_j^k)^*} \rukl{\ilim_j \gekl{C_0(X_{k,j}); p_{k+j}^*}} = C_0(\bX_k).
\end{align}

For all $k$, define ${\bar \imath}_k: \: \bSigma_k \to \bSigma_{k+1}, \, (\sigma_j)_j \ma (\imath_{k+j}(\sigma_{j+1}))_j$. Clearly, ${\bar \imath}_k$ is an injective twisted groupoid homomorphism. The image of ${\bar \imath}_k$ is given by ${\bar \imath}_k (\bSigma_k) = (p_0^{k+1})^{-1}(T_k)$, hence it is open. Similarly, ${\bar i}_k: \: \bG_k \to \bG_{k+1}, \, (\gamma_j)_j \ma (i_{k+j}(\gamma_{j+1}))_j$ is an injective continuous groupoid homomorphism with open image. Let $\bSigma \defeq \ilim_k \gekl{\bSigma_k; {\bar \imath}_k}$ and $\bG \defeq \ilim_k \gekl{\bG_k; {\bar i}_k}$, equipped with the inductive limit topology. Let ${\bar \imath}_k^{\infty}: \: \bSigma_k \to \bSigma$ and ${\bar i}_k^{\infty}: \: \bG_k \to \bG$ be the canonical inclusions. The maps ${\bar \imath}_k^{\infty}$ have open image, as $({\bar \imath}_l^{\infty})^{-1}(\bSigma_k)  = {\bar \imath}_{l,k}(\bSigma_k) \subseteq \bSigma_l$ is open for $l \geq k$, where ${\bar \imath}_{l,k} = {\bar \imath}_{l-1} \circ \dotso \circ {\bar \imath}_{k+1} \circ {\bar \imath}_k$. Similarly, the ${\bar i}_k^{\infty}$ all have open image.

By construction, it is clear that $(\bG,\bSigma)$ is a twisted \etale{} groupoid. It is also clear that $\bG^{(0)} = \ilim_k \gekl{\bX_k; {\bar \imath}_k}$. By Lemma~\ref{Lem:i:opensubgpd-hom}, ${\bar \imath}_k^{\infty}$ induces an isometric homomorphism $({\bar \imath}_k^{\infty})_*: \: C^*_r(\bG_k,\bSigma_k) \to C^*_r(\bG,\bSigma)$. We hence obtain an isometric homomorphism $\ilim_k ({\bar \imath}_k^{\infty})_*: \: \ilim_k \gekl{C^*_r(\bG_k,\bSigma_k); ({\bar \imath}_k)_*} \to C^*_r(\bG,\bSigma)$. It has dense image since every $f \in C_c(\bG,\bSigma)$ has $\supp(f) \subseteq {\bar \imath}_k^{\infty}(\bG_k)$ for some $k$, so that $f \in ({\bar \imath}_k^{\infty})_*(C^*_r(\bG_k,\bSigma_k))$. Hence $\ilim_k ({\bar \imath}_k^{\infty})_*$ is an isomorphism $\ilim_k \gekl{C^*_r(\bG_k,\bSigma_k); ({\bar \imath}_k)_*} \overset{\cong}{\lori} C^*_r(\bG,\bSigma)$. Since $\bG^{(0)} = \ilim_k \gekl{\bX_k; {\bar \imath}_k}$, we must have
\begin{align}
\label{ilimbX_k}
  {\ilim_k ({\bar \imath}_k^{\infty})_*} \rukl{\ilim_k \gekl{C_0(\bX_k); ({\bar \imath}_k)_*}} = C_0(\bG^{(0)}).
\end{align}

Now consider the composition $({\bar \imath}_k^{\infty})_* \circ (p_0^k)^*: \: C^*_r(G_k,\Sigma_k) \overset{(p_0^k)^*}{\lori} C^*_r(\bG_k,\bSigma_k) \overset{({\bar \imath}_k^{\infty})_*}{\lori} C^*_r(\bG,\bSigma)$. As
$$
  ({\bar \imath}_{k+1}^{\infty})_* \circ (p_0^{k+1})^* \circ (\imath_k)_* \circ (p_k)^* = ({\bar \imath}_k^{\infty})_* \circ (p_0^k)^*,
$$
we obtain an isometric homomorphism
$$
  \ilim_k \rukl{({\bar \imath}_k^{\infty})_* \circ (p_0^k)^*}: \: \ilim_k \gekl{C^*_r(G_k,\Sigma_k); (\imath_k)_* \circ (p_k)^*} \to C^*_r(\bG,\bSigma).
$$
It has dense image because it contains, for all $j$,
\begin{align*}
  \rukl{({\bar \imath}_{k+j}^{\infty})_* \circ (p_0^{k+j})^*} \rukl{C_c(G_{k,j},\Sigma_{k,j})}
  &= \rukl{({\bar \imath}_{k+j}^{\infty})_* \circ ({\bar \imath}_{k+j,k})_* \circ (p_j^k)^*} \rukl{C_c(G_{k,j},\Sigma_{k,j})}\\
  &= \rukl{({\bar \imath}_k^{\infty})_* \circ (p_j^k)^*} \rukl{C_c(G_{k,j},\Sigma_{k,j})},
\end{align*}
so that it is dense in $({\bar \imath}_k^{\infty})_* (C^*_r(\bG_k,\bSigma_k))$ for all $k$, hence also dense in $C^*_r(\bG,\bSigma)$.

Thus $\ilim_k \rukl{({\bar \imath}_k^{\infty})_* \circ (p_0^k)^*}$ is an isomorphism $\ilim_k \gekl{C^*_r(G_k,\Sigma_k); (\imath_k)_* \circ (p_k)^*} \overset{\cong}{\lori} C^*_r(\bG,\bSigma)$. \eqref{ilimX_kj} and \eqref{ilimbX_k} show that $\rukl{\ilim_k \rukl{({\bar \imath}_k^{\infty})_* \circ (p_0^k)^*}} \rukl{\ilim_k \gekl{C_0(X_k); (\imath_k)_* \circ (p_k)^*}} = C_0(\bG^{(0)})$. This proves (a) and (b).

For (c), assume that for every $k$, $G_k$ is topologically principal. We first show that
$$
  \fX_k \defeq \menge{{\bar x} = (x_j) \in \bX_k}{(G_{k+i})^{x_i}_{x_i} = \gekl{x_i} \ {\rm for} \ {\rm some} \ i}
$$
is dense in $\bX_k$. As non-empty open subsets of the form $\bigcap_{h=0}^i (p_h^k)^{-1}(U_h)$ form a basis for the topology of $\bX_k$, where the $U_h$ are open subsets of $X_{k+h}$, it suffices to find ${\bar x} \in \fX_k$ with ${\bar x} \in \bigcap_{h=0}^i (p_h^k)^{-1}(U_h)$. We have $\bigcap_{h=0}^i (p_h^k)^{-1}(U_h) = (p_i^k)^{-1} \rukl{\bigcap_{h=0}^i p_{h,i}^{-1}(U_h)}$, where $p_{h,i} = p_{k+h} \circ \dotso \circ p_{k+i-1}$. As $G_{k+i}$ is topologically principal, we can find $x_i \in \bigcap_{h=0}^i p_{h,i}^{-1}(U_h)$ with $(G_{k+i})^{x_i}_{x_i} = \gekl{x_i}$. As $p_i^k$ is surjective, we can find ${\bar x} \in \bX_k$ with $p_i^k({\bar x}) = x_i$. Thus, $\fX_k$ is indeed dense in $\bX_k$. Hence, for all $k$, ${\bar \imath}_k^{\infty}(\fX_k)$ is dense in ${\bar \imath}_k^{\infty}(\bX_k)$. It thus remains to show that for every $k$ and ${\bar x} \in \fX_k$, ${\bar \imath}_k^{\infty}({\bar x})$ has trivial stabilizer group in $\bG$. As ${\bar x} = (x_j)$ lies in $\fX_k$, there exists $i$ such that $(G_{k+i})^{x_i}_{x_i} = \gekl{x_i}$. Let us first show that for all $j \geq i$, we must have $(G_{k+j})^{x_j}_{x_j} = \gekl{x_j}$. Proceeding inductively on $j$, we may assume that $(G_{k+j})^{x_j}_{x_j} = \gekl{x_j}$ and have to show that $(G_{k+j+1})^{x_{j+1}}_{x_{j+1}} = \gekl{x_{j+1}}$. By assumption, $\dot{p}_{k+j}$ restricts to a bijection $(G_{k+j+1})_{x_{j+1}} \to (G_{k+j})_{x_j}$. Take $\gamma \in (G_{k+j+1})_{x_{j+1}}$ with $r(\gamma) = x_{j+1}$. Then $r(\dot{p}_{k+j}(\gamma)) = x_j$, so that $\dot{p}_{k+j}(\gamma) = x_j$, and thus $\gamma = x_{j+1}$. Hence, we indeed have $(G_{k+j})^{x_j}_{x_j} = \gekl{x_j}$ for all $j \geq i$. Now take ${\bar \gamma} \in \bG$ with range and source equal to ${\bar \imath}_k^{\infty}({\bar x})$. We have ${\bar \gamma} = {\bar i}_l^{\infty}(\gamma)$ for some $\gamma \in \bG_l$, and we may without loss of generality assume that $l \geq k$. Let $y = {\bar \imath}_{l,k}({\bar x})$. It follows from injectivity of ${\bar i}_l^{\infty}$ that $\gamma \in (\bG_l)^y_y$. Write $y = (y_j)$ and $\gamma = (\gamma_j)$. Then we have $y_j = x_{j+l-k}$ for all $j$. Hence, by what we have just shown, $(G_{l+j})_{y_j}^{y_j} = \gekl{y_j}$ for all $j \geq i$. This implies $\gamma_j = y_j$ for all $j \geq i$. But then we must have $\gamma_j = y_j$ for all $j$, since $\gamma_h$ for $h \leq i$ is determined recursively by $\gamma_{h-1} = \dot{p}_{l+h-1}(\gamma_h)$. Hence $\gamma = y$, and we conclude that ${\bar \gamma} = {\bar \imath}_k^{\infty}({\bar x})$. Thus $\bG$ is indeed topologically principal. This shows (c).
\eproof

\bremark
\label{Rem:trivialSigmaIndLim}
If the twists $\Sigma_k$ are all trivial (i.e., $\Sigma_k = \Tz \times G_k$) in Theorem~\ref{indlim--GPD}, then also $\bSigma_k$ and $\bSigma$ are trivial. This is easy to see from the proof.
\eremark

\bremark
Let us explain how \cite[Lemma~5.2]{BL} may be derived as a special case of Theorem~\ref{indlim--GPD}. Recall the setting in \cite[Lemma~5.2]{BL}: For all $j = 1, 2, \dotsc$, $(G(j),\Sigma(j))$ is a twisted \etale{} topologically principal groupoid with compact unit space $X(j)$. Let
$$A_k \defeq C^*_r(G(1),\Sigma(1)) \otimes_{\min} \dotso \otimes_{\min} C^*_r(G(k),\Sigma(k)).$$
Consider $\vp_k: \: A_k \to A_{k+1}, \, a \ma a \otimes 1$. Let $B_k \defeq C(X(1)) \otimes \dotso \otimes C(X(k)) \subseteq A_k$. We want to see that $\rukl{\ilim_k \gekl{A_k; \vp_k}, \ilim_k \gekl{B_k; \vp_k}}$ is a Cartan pair. In order to apply Theorem~\ref{indlim--GPD}, set $G_k \defeq G(1) \times \dotso \times G(k)$, $X_k \defeq X(1) \times \dotso \times X(k)$ and $\Sigma_k \defeq \Sigma(1) \times_{\Tz} \dotso \times_{\Tz} \Sigma(k)$, where $\Sigma(1) \times_{\Tz} \Sigma(2)$ denotes the quotient of $\Sigma(1) \times \Sigma(2)$ with respect to $(z \sigma,\sigma') \sim (\sigma, z \sigma')$ (as in \cite[\S~5]{BL}). It is easy to see that $(G_k,\Sigma_k)$ is a twisted groupoid such that $(A_k,B_k) \cong (C^*_r(G_k,\Sigma_k),C(X_k))$. Moreover, set $H_k \defeq G(1) \times \dotso \times G(k) \times X(k+1)$. Obviously, $H_k$ is an open subgroupoid in $G_{k+1}$. Similarly, $T_k \defeq \Sigma(1) \times \dotso \times \Sigma(k) \times X(k+1)$ is an open subgroupoid of $\Sigma_{k+1}$. Let $\imath_k: \: T_k \into \Sigma_{k+1}$ and $i_k: \: H_k \into G_{k+1}$ be the canonical inclusions. Define $p_k: \: T_k \to \Sigma_k, \, (\vec{\sigma},x) \ma \vec{\sigma}$ for $\vec{\sigma} \in \Sigma_k$. It is easy to check that $p_k$ is a surjective and proper twisted groupoid homomorphism, for all $k$. Clearly, each restriction $(H_k)_x  \to (G_k)_{p_k(x)}$, $\eta \mapsto \dot{p}_k(\eta)$ is a bijection. Moreover, we have $\vp_k = (i_k)_* \circ (p_k)^*$. As conditions (LT) and (E) are easily verified, Theorem~\ref{indlim--GPD} applies and yields that
$$
  \rukl{\ilim_k \gekl{A_k; \vp_k}, \ilim_k \gekl{B_k; \vp_k}}
  \cong \rukl{\bigotimes_{j=1}^{\infty} C^*_r(G(j),\Sigma(j)), \bigotimes_{j=1}^{\infty} C(X(j))} \cong C^*_r(\bG,\bSigma)
$$
is a Cartan pair.
\eremark

\section{Invariant Cartan subalgebras and the UCT problem}
\label{Sec:Inv}

Recall that an inverse semigroup is a semigroup $\cS$ with the property that for every $s \in \cS$, there exists a unique $s^* \in \cS$ such that $s = ss^*s$ and $s^* = s^*ss^*$. We write
\[
E := \menge{e \in \cS}{e = e^2} = \menge{e \in \cS}{e = e^* = e^2}
\]
for the associated semilattice of idempotents. If $A$ is a C*-algebra and $\cS \subseteq A$ is an inverse semigroup of partial isometries in $A$, then we write $C^*(E)$ for the sub-C*-algebra of $A$ generated by $E$. We remark that in the literature $C^*(E)$ often denotes the universal C*-algebra of the semilattice $E$. However, there should be no confusion as the latter does not appear anywhere in this paper. 

\blemma \label{Spec(B) tot disc}
Let $(A,B)$ be a Cartan pair with $A$ a separable C*-algebra and such that $\Spec(B)$ (i.e., the space of characters of the commutative C*-algebra $B$) is totally disconnected. Then $\cS = \menge{s \in N_A(B)}{s=ss^*s}$ is an inverse semigroup generating $A$ as a C*-algebra such that $B = C^*(E)$. Moreover, if $\alpha \in \Aut(A)$ is an automorphism, then $\alpha(B) = B$ if and only if $\alpha(\cS) = \cS$.
\elemma
\bproof
It is clear from the definition of the normalizer that $\cS$ is an inverse semigroup. For instance, to see that $\cS$ is multiplicatively closed, take $s, t \in \cS$. Then $s$ and $t$ are partial isometries with commuting source and range projections as these projections lie in $B$ because $s$ and $t$ lie in $N_A(B)$. It follows that $(st)(st)^*$ is idempotent because $((st)(st)^*)^2 = s t t^* s^* s t t^* s^* = s (t t^*) (s^* s) (t t^*) s^* = s (t t^*)^2 s^* s s^* = (st)(st)^*$, and thus $st$ is again a partial isometry. Since we know that $st$ lies in $N_A(B)$, it follows that $st$ lies in $\cS$. Let us prove that $\cS$ generates $A$ as a C*-algebra. Let $(G,\Sigma)$ be a twisted \etale{} Hausdorff locally compact second countable topologically principal groupoid with $(A,B) \cong (C^*_r(G,\Sigma),C_0(G^{(0)}))$, which exists by Theorem~\ref{Renault}. We claim that it is enough to prove that any element in $N_{C^*_r(G,\Sigma)}(C_0(G^{(0)}))$ with compact open support can be approximated by linear combinations of partial isometries in $N_{C^*_r(G,\Sigma)}(C_0(G^{(0)}))$. To see this, observe that it is enough to consider elements in $N_{C^*_r(G,\Sigma)}(C_0(G^{(0)}))$ with compact support, as these are dense in $N_{C^*_r(G,\Sigma)}(C_0(G^{(0)}))$ by \cite[Lemma~5.8]{R}. By \cite[Proposition~4.8]{R}, $N_{C^*_r(G,\Sigma)}(C_0(G^{(0)}))$ consists exactly of the sections whose open support is a bisection. As $G$ is \etale{} and $G^{(0)}$ is totally disconnected, the compact open bisections of $G$ form a basis for the topology. Given $f \in N_{C^*_r(G,\Sigma)}(C_0(G^{(0)}))$ with compact support, we can therefore find a finite covering of $\supp(f)$ by finitely many compact open bisections. Using a subordinate partition of unity, we can thus write $f$ as a linear combination of elements in $N_{C^*_r(G,\Sigma)}(C_0(G^{(0)}))$ with compact open support. 

Assume therefore that $f$ has compact open support $K \subseteq G$ and norm one. By \cite[Proposition~4.8]{R}, $s(K) = \dom(f) = \menge{x \in G^{(0)}}{f^*f(x) > 0}$, which is a compact open subset of $G^{(0)}$. Let $\epsilon > 0$ and find pairwise disjoint compact open subsets $U_1,\ldots,U_n$ of $G^{(0)}$ whose union is $s(K)$, such that
$$
  \abs{\frac{f^*f(x) - f^*f(y)}{f^*f(x)}} < \epsilon
$$
for all $x,y \in U_j$ and $j=1,\ldots,n$. For each $j$, choose a fixed $x_j \in U_j$, and define $g_j \in C_0(G^{(0)})$ with $\supp(g_j) \subseteq U_j$ by setting
$$
g_j(x) \defeq \frac{f^*f(x_j)}{f^*f(x)}.
$$
Then $\abs{1 - g_j(x)} < \epsilon$ for all $x \in U_j$ and $f^*fg_j$ is constant on $U_j$. So $fg_j^{1/2} \in N_{C^*_r(G,\Sigma)}(C_0(G^{(0)}))$ is a scalar multiple of a partial isometry. Now, $\sum_{j=1}^n fg_j^{1/2}$ is a linear combination of partial isometries in $N_{C^*_r(G,\Sigma)}(C_0(G^{(0)}))$ satisfying 
\[
\| f - \sum_{j=1}^n fg_j^{1/2}\| \leq \| \chi_{s(K)} - \sum_{j=1}^n g_j^{1/2} \|_\infty < \epsilon.
\]

By definition of $\cS$, it holds that $E$ contains all projections of $B$. As $B$ has totally disconnected spectrum, we conclude that $B = C^*(E)$.

Now let $\alpha \in \Aut(A)$. If $\alpha(\cS) = \cS$, then $\alpha(E) = E$ and consequently $\alpha(B) = B$. Conversely, if $\alpha(B) = B$, then always $\alpha(N_A(B)) \subseteq N_A(B)$ and consequently $\alpha(\cS) \subseteq \cS$. As $\alpha^{-1}(B) = B$ as well, it follows that $\alpha(\cS) = \cS$.
This concludes the proof.
\eproof

\bremark
\label{Rem1:tight}
If $A$ in the situation of Lemma~\ref{Spec(B) tot disc} is unital and nuclear, and the Weyl twist of $(A,B)$ (see Remark~\ref{A,B-->GPD}) is trivial, then $A$ is canonically isomorphic to the tight C*-algebra of some inverse semigroup in the sense of \cite[Definition~2.2]{Ex}. For this, consider the inverse subsemigroup $\cS' \subseteq \cS$ consisting of those partial isometries in $N_A(B)$, which correspond to indicator functions on compact open bisections of $G$. By \cite[Theorem~4.9]{Ex}, the canonical inclusion $\cS' \into A$ induces an isomorphism $A \cong C^*_{tight}(\cS')$. In fact, this remains true if we replace $\cS'$ by a possibly smaller inverse subsemigroup, as long as the corresponding compact open bisections form a basis for the topology of $G$.
\eremark

The following result is a slight improvement of \cite[Theorem~2.3 with proof]{BS}. We are grateful to G{\'a}bor Szab{\'o} for pointing out to us that only the separability assumption on the C*-algebra in the statement of \cite[Theorem~2.3]{BS} is needed. Given two automorphisms $\alpha$ and $\beta$ of a C*-algebra $A$, we write $\alpha \ue \beta$ if $\alpha$ and $\beta$ are approximately unitarily equivalent via unitaries in the minimal unitization of $A$.

\btheo \label{thm:result BS}
Let $\Gamma$ be a finite group. Let $A$ be a separable C*-algebra such that $A \cong M_{|\Gamma|^\infty} \otimes A$. Let $\menge{\beta_\gamma}{\gamma \in \Gamma}$ be a family of *-automorphisms on $A$ such that $\beta_{\gamma \eta} \ue \beta_\gamma \circ \beta_\eta$ for all $\gamma,\eta \in \Gamma$. Consider the C*-dynamical system
\[
(D,\delta) := \lim\limits_{k \to \infty} \gekl{(C(\Gamma) \otimes M^{\otimes k-1} \otimes A, \sigma \otimes \id_{M^{\otimes k-1} \otimes A}), \vp_k},
\]
where $M$ denotes the algebra of $|\Gamma|\times |\Gamma|$ matrices, $\sigma:\Gamma \curvearrowright C(\Gamma)$ is the action given by $\sigma_\gamma(f)(\eta) = f(\gamma^{-1}\eta)$ and
\[
\vp_k(f)(\gamma) := \sum\limits_{\eta \in \Gamma} e_{\eta,\eta} \otimes (\id \otimes \beta_{\eta})(f(\gamma\eta)),
\]
under the canonical identifications $C(\Gamma) \otimes M^{\otimes r} \otimes A \cong C(\Gamma, M^{\otimes r} \otimes A)$. Then the action $\delta:\Gamma \curvearrowright D$ has the Rokhlin property and there exists a *-isomorphism $\mu: D \to A$ such that $\mu \circ \delta_\gamma \circ \mu^{-1} \ue \beta_\gamma$ for all $\gamma \in \Gamma$.
\etheo
\bproof
If we assume that $A$ is either unital, stable or has stable rank one, then this follows from \cite[Theorem~2.3 and its proof]{BS}. The only reason for this restriction is that in these cases the two notions of approximate unitary equivalence, with respect to unitaries in the multiplier algebra and the minimal unitization, coincide. However, apart from unitaries in $C(\Gamma) \otimes M^{\otimes k} \otimes 1 \subseteq \cM(C(\Gamma) \otimes M^{\otimes k} \otimes A)$, which are used in the proof of \cite[Theorem~2.3]{BS} (in the fourth paragraph, for $\Psi_n \approx_u {\rm Ad} \rukl{\bigoplus_{\gamma \in \Gamma} \lambda(\gamma^{-1}) \otimes \mathbf{1}_{A^{(n)}}} \circ \Psi_n$), the only point where multiplier unitaries are used is in \cite[Remark~1.7]{BS}. All these unitaries can be chosen to be norm-homotopic to the identity. However, the connected component of the identity in the unitary group of the minimal unitization of a C*-algebra $C$ is always strictly dense in the connected component of the unitary group of its multiplier algebra $\cM(C)$. This is a consequence of the fact that the set of self-adjoint element in $C$ is stictly dense in the set of self-adjoint elements of $\cM(C)$. This shows that all unitaries in the proof of \cite[Theorem~2.3]{BS} implementing approximate unitary equivalences can be taken in the respective minimal unitizations.
\eproof

For a discrete group action $\alpha:\Gamma \curvearrowright A$ on a C*-algebra, we write $t_\gamma \in \cM(A \rtimes_\alpha \Gamma)$ for the canonical unitary implementing $\alpha_\gamma$.

We say that an element $a \in A$ is $\alpha$-homogeneous if for every $\gamma \in \Gamma$, there is some scalar $\lambda \in \Cz$ such that $\alpha_\gamma(a) = \lambda a$. Note that $\lambda$ is of modulus one if $a \neq 0$, with $\lambda^n = 1$ if $\alpha_\gamma^n = \id$. Furthermore, a subset $X \subseteq A$ is said to be $\alpha$-homogeneous if it consists of $\alpha$-homogeneous elements.

\btheo
\label{Cartan fixed general}
Let $\Gamma$ be a finite group. Let $A$ be a separable C*-algebra with $A \cong M_{|\Gamma|^\infty} \otimes A$. Let $\alpha:\Gamma \curvearrowright A$ be an action with the Rokhlin property. Assume that there exists a family of *-automorphisms $\menge{\beta_\gamma}{\gamma \in \Gamma} \subseteq \Aut(A)$ and a Cartan subalgebra $B \subseteq A$ such that $\beta_\gamma \ue \alpha_\gamma$ and $\beta_\gamma(B) = B$ for all $\gamma \in \Gamma$.

Then there exists another Cartan subalgebra $C \subseteq A $ such that
\begin{itemize}
\item[1)] $\alpha_\gamma(C)=C$ for all $\gamma \in \Gamma$;
\item[2)] $C \subseteq A \rtimes_\alpha \Gamma$ is also a Cartan subalgebra, which is fixed pointwise by $\hat{\alpha}$ if $\Gamma$ abelian;
\item[3)] if $A$ is unital and $\Gamma$ abelian, then unitaries in $A = (A \rtimes_\alpha \Gamma)^{\hat{\alpha}}$ witnessing approximate representability of $\hat{\alpha}$ can be chosen in $C$;
\item[4)] if $B$ has totally disconnected spectrum, then there exists an $\alpha$-invariant inverse semigroup $\cS \subseteq A$ of partial isometries generating $A$ as a C*-algebra such that $C^*(E) = C$.
In this case, $\tilde{\cS} \defeq \menge{st_\gamma}{\gamma \in \Gamma,\ s \in \cS}$ is an inverse semigroup of partial isometries in $A \rtimes_{\alpha} \Gamma$ which generates $A \rtimes_\alpha \Gamma$ as a C*-algebra and whose idempotent semilattice is equal to $E$.
If $\Gamma$ is abelian, then $\tilde{\cS}$ consists of $\hat{\alpha}$-homogeneous elements;
\item[5)] if $A$ is unital and the canonical map $\cP(B)\setminus \gekl{0,1} \to K_0(A)$ is surjective, then so is the canonical map $\cP(C) \setminus \gekl{0,1} \to K_0(A)$.
\end{itemize}

If $A$ is either a unital UCT Kirchberg algebra or a unital, simple, separable, nuclear TAF algebra in the UCT class, then the conclusions remain true if we replace $\beta_\gamma \ue \alpha_\gamma$ by $K_*(\beta_\gamma)=K_*(\alpha_\gamma)$ for all $\gamma \in \Gamma$.
\etheo
\bproof
By Theorem~\ref{thm:result BS}, the action $\delta:\Gamma \curvearrowright D$ given by
\[
(D,\delta) := \lim\limits_{k \to \infty} \gekl{(A_k, \sigma \otimes \id_{M^{\otimes k-1} \otimes A}), \vp_k},
\]
where $A_k = C(\Gamma) \otimes M^{\otimes k-1} \otimes A$ and $\vp_k: \: A_k \to A_{k+1}$ satisfies
$$
  \vp_k(f)(\gamma) = \sum_{\eta \in \Gamma} e_{\eta,\eta} \otimes (\id \otimes \beta_{\eta})(f(\gamma \eta)),
$$
has the Rokhlin property and there exists an isomorphism $\mu:D \to A$ such that $\mu \circ \delta_\gamma \circ \mu^{-1} \ue \beta_\gamma$ for all $\gamma \in \Gamma$. By assumption, this implies that $\mu \circ \delta_\gamma \circ \mu^{-1} \ue \alpha_\gamma$ for all $\gamma \in \Gamma$. Hence, $(A,\alpha)$ and $(D,\delta)$ are conjugate by \cite[Theorem~3.5 and Remark~3.6]{Na} (see also \cite[Theorem~3.5]{Iz}). We show that $D$ has a Cartan subalgebra satisfying 1) - 5) for $\delta$.

By Theorem~\ref{Renault}, there exists a twisted \etale{} topologically principal groupoid $(G,\Sigma)$ such that $(A,B) \cong (C^*_r(G,\Sigma),C_0(G^{(0)}))$. Denote by $X \defeq G^{(0)}$ the unit space and by $\pi: \: \Sigma \twoheadrightarrow G$ the canonical projection. As $\beta_\gamma(B) = B$ for given $\gamma \in \Gamma$, there exists a twisted groupoid automorphism $q_{\gamma}: \: \Sigma \overset{\cong}{\lori} \Sigma$ such that $\beta_\gamma = (q_{\gamma})^*$ under the above identification, see \cite[proof of Proposition~3.4]{BL}. Let $R = \Gamma \times \Gamma$, viewed as a relation, i.e., a groupoid, so that $C^*_r(R) = M$. The unit space of $R$ is given by the diagonal $\Delta$, which is canonically identified with $\Gamma$. Consider the groupoid $G_k \defeq \Gamma \times R^{k-1} \times G$. The first copy of $\Gamma$ is viewed as a disjoint union of points (i.e., one-point groupoids). $G_k$ becomes an \etale{} topologically principal groupoid under component-wise operations, and its unit space is $X_k = \Gamma \times \Delta^{k-1} \times X$. Set $\Sigma_k \defeq \Gamma \times R^{k-1} \times \Sigma$. With respect to the projection $\pi_k \defeq \id_{\Gamma \times R^{k-1}} \times \pi$, $(G_k,\Sigma_k)$ is a twisted \etale{} topologically principal groupoid. We have a canonical identification $A_k \cong C^*_r(G_k,\Sigma_k)$ which identifies $B_k \defeq C(\Gamma) \otimes \tilde{D}^{\otimes (k-1)} \otimes B$ with $C_0(X_k)$. Here $\tilde{D}$ is the canonical diagonal in $M$. Furthermore, set $H_k \defeq \Gamma \times \Delta \times R^{k-1} \times G$. Obviously, $H_k$ is an open subgroupoid in $G_{k+1}$, with unit space $Y_k \defeq \Gamma \times \Delta \times \Delta^{k-1} \times X$. Let $i_k: \: H_k \into G_{k+1}$ be the canonical inclusion. Set $T_k \defeq \Gamma \times \Delta \times R^{k-1} \times \Sigma$. $(H_k,T_k)$ is a twisted \etale{} topologically principal groupoid. $T_k$ is an open subgroupoid of $\Sigma_{k+1}$. Let $\imath_k: \: T_k \into \Sigma_{k+1}$ be the canonical inclusion. Define $p_k: \: T_k \to \Sigma_k, \, (\gamma, d, \bfr, \sigma) \ma (\gamma d, \bfr, q_d(\sigma))$. Clearly, $p_k$ is a twisted groupoid homomorphism, which is surjective, proper, and $\dot{p}_k$ restricts to bijections $(H_k)_v \cong (G_k)_{p_k(v)}$. Moreover, it is easy to check that $\vp_k = (i_k)_* \circ (p_k)^*$. Furthermore, it is easy to see that conditions (LT) and (E) are satisfied. Thus Theorem~\ref{indlim--GPD} applies and yields that
$$\rukl{\ilim_k \gekl{A_k; \vp_k}, \ilim_k \gekl{B_k; \vp_k}} \cong \rukl{C^*_r(\bG, \bSigma), C_0(\bG^{(0)})}$$
is a Cartan pair. As $\sigma \otimes \id_{M^{\otimes k-1} \otimes A}(B_k) = B_k$ for all $k \geq 1$, the Cartan subalgebra $F := \ilim_k \gekl{B_k; \vp_k}$ of $D$ gets fixed globally by $\delta$. This shows 1).

Set $\cG \defeq \Gamma \rtimes \Gamma$, where we let $\Gamma$ act on itself by left multiplication. Set $G_k \defeq \cG \times R^{k-1} \times G$, $\Sigma_k \defeq \cG \times R^{k-1} \times \Sigma$, $\pi_k \defeq \id_{\cG \times R^{k-1}} \times \pi$,  $H_k \defeq \cG \times \Delta \times R^{k-1} \times G$, $T_k \defeq \cG \times \Delta \times R^{k-1} \times \Sigma$, and $p_k: \: T_k \to \Sigma_k, \, (\bfg, d, \bfr, \sigma) \ma (\bfg \cdot d, \bfr, \sigma)$, where for $\bfg = (\alpha,\gamma) \in \cG$, $\bfg \cdot d = (\alpha d,\gamma)$. Applying Theorem~\ref{indlim--GPD} to this setting gives us that
$$(D \rtimes_\delta \Gamma,F) = \rukl{\ilim_k \gekl{A_k \rtimes \Gamma; \vp_k \rtimes \Gamma}, \ilim_k \gekl{B_k; \vp_k \rtimes \Gamma}}$$
is a Cartan pair. We conclude 2).

Let $A$ be unital and $\Gamma$ abelian. The above construction of the Cartan pair ($D$,$F$) shows that the Rokhlin projections for $\delta$ can all be chosen to be in the subalgebras $C(\Gamma) \subseteq B_k$, $k \geq 1$. In particular, all Rokhlin projections for $\delta$ can be chosen in $F$. It now follows from Remark~\ref{rem: unitaries approx repr} that the unitaries in $D = (D \rtimes_\delta \Gamma)^{\hat{\delta}}$  witnessing approximate representability of $\hat{\delta}$ can also be chosen to be in $F$. This shows 3).

Assume now that $\Spec(B)$ is totally disconnected. Then $B$ is generated by its projections. The construction of the Cartan pair $(D,F)$ shows that in this case $F$ is generated by its projections as well and therefore has totally disconnected spectrum.
Applying Lemma~\ref{Spec(B) tot disc}, we obtain an $\alpha$-invariant inverse semigroup $\cS \subseteq D$ of partial isometries generating $D$ as a C*-algebra such that $F = C^*(E)$. The property that $\delta$ preserves $F$ is equivalent to the canonical unitaries $t_\gamma$, $\gamma \in \Gamma$, normalizing $F$. Therefore, $\tilde{\cS} = \menge{st_\gamma}{\gamma \in \Gamma,\ s \in \cS}$ is an inverse subsemigroup of partial isometries in $N_{D \rtimes_\delta \Gamma}(F)$ with idempotent semilattice equal to $E$ such that $(C^*(\tilde{\cS}),C^*(E)) = (D \rtimes_{\delta} \Gamma,F)$.
If $\Gamma$ is abelian, then $\tilde{\cS}$ clearly consists of $\hat{\delta}$-homogeneous elements. This shows 4).

Assume that $A$ is unital and $\cP(B) \setminus \gekl{0,1} \to K_0(A)$ is surjective. Let $x \in K_0(D)$ and find some $k \in \Nz$ and $\tilde{x} \in K_0(C(\Gamma) \otimes M^{\otimes k} \otimes A)$ which gets mapped to $x$ under the canonical map. Using the identification
\[
K_0(C(\Gamma) \otimes M^{\otimes k} \otimes A) \cong K_0(A)^{|\Gamma|},
\]
we find by assumption $q_1,\ldots,q_{|\Gamma|} \in \cP(B) \setminus \gekl{0,1}$ such that $\tilde{x} = ([q_1],\ldots,[q_{|\Gamma|}])$. This shows that $x$ is in the image of the canonical map $\cP(C(\Gamma) \otimes M^{\otimes k} \otimes B) \setminus \gekl{0,1} \to K_0(D)$. Using that the canonical map $C(\Gamma) \otimes M^{\otimes k} \otimes A \to D$ is injective, it follows that $\cP(F) \setminus \gekl{0,1} \to K_0(D)$ is surjective. This shows 5).

Finally, assume that $A$ is either a unital UCT Kirchberg algebra or a unital, simple, separable, nuclear TAF algebra in the UCT class. Assume that the family $\menge{\beta_\gamma}{\gamma \in \Gamma} \subseteq \Aut(A)$ satisfies $K_*(\beta_\gamma) = K_*(\alpha_\gamma)$ and $\beta_\gamma(B) = B $ for all $\gamma \in \Gamma$. By Theorem~\ref{thm:result BS}, we find an action with the Rokhlin property $\rho:\Gamma \curvearrowright A$ such that $\beta_\gamma \ue \rho_\gamma$ for all $\gamma \in \Gamma$. In particular, $K_*(\alpha_\gamma) = K_*(\rho_\gamma)$ for all $\gamma \in \Gamma$. By the first part of this theorem, we obtain a Cartan subalgebra $F \subseteq A$ satisfying 1) - 5) (if applicable) for $\rho$. By \cite[Theorem~4.2 and Theorem~4.3]{Iz2}, there exists an equivariant isomorphism $(A,\alpha) \cong (A,\rho)$. The image $C$ of $F$ under this equivariant isomorphism yields the desired Cartan subalgebra.
\eproof

\bremark
\label{Rem2:tight}
If the Weyl twist of $(A,B)$ in Theorem~\ref{Cartan fixed general} is trivial, then also the Weyl twists of $(A,C)$ and $(A \rtimes_{\alpha} \Gamma,C)$ are trivial. This follows from the above proof and Remark~\ref{Rem:trivialSigmaIndLim}. If in addition $A$ is unital and nuclear, it then follows from Remark~\ref{Rem1:tight} that for $\cS \subseteq A$ and $\tilde{\cS} \subseteq A \rtimes_{\alpha} \Gamma$ in Theorem~\ref{Cartan fixed general}~4), there exist inverse subsemigroups $\cS' \subseteq \cS$ and $\tilde{\cS}' \subseteq \tilde{\cS}$ such that the inclusions $\cS' \into A$ and $\tilde{\cS}' \into A \rtimes_{\alpha} \Gamma$ induce isomorphisms $A \cong C^*_{tight}(\cS')$ and $A \rtimes_{\alpha} \Gamma \cong C^*_{tight}(\tilde{\cS}')$. As the action $\alpha$ is induced by a twisted groupoid automorphism, $\cS'$ is $\alpha$-invariant; see Lemma~\ref{Spec(B) tot disc}.
\eremark

Let us now explain how the following result is a consequence of \cite{Kat08a,Sp07a,Sp07b}. Recall that a subgroup of a group is called a Sylow subgroup if it is a maximal $p$-subgroup for some prime number $p \geq 2$. Given a countable row-finite graph $E$, we denote by $\cO(E)$ the graph C*-algebra associated with $E$.
\blemma
\label{Lem:Sp+Kat}
Let $K_0$, $K_1$ be countable abelian groups and $u \in K_0$ an arbitrary element. Let $\Gamma$ be a finite group such that every Sylow subgroup of $\Gamma$ is cyclic. Let $\Gamma \curvearrowright K_i$, $i=0,1$, be actions such that $\Gamma \curvearrowright K_0$ fixes $u$.

Then there exists an \etale{} topologically principal groupoid $G$ with compact totally disconnected unit space $X$ and an action $\Gamma \curvearrowright G$ such that
\begin{enumerate}
\item $C^*_r(G)$ is a unital UCT Kirchberg algebra.
\item We have a $\Gamma$-equivariant isomorphism $(K_0(C^*_r(G)),[1],K_1(C^*_r(G))) \cong (K_0,u,K_1)$.
\item The canonical map $\cP(C(X)) \setminus \gekl{0,1} \to K_0(C^*_r(G))$ is surjective.
\end{enumerate}
\elemma
\bproof
Our proof is a combination of \cite{Kat08a,Sp07a,Sp07b}. We first use \cite{Kat08a} to obtain exact sequences of permutation modules. We then apply the construction in \cite{Sp07b}, which builds on \cite{Sp07a}. This approach has also been mentioned in \cite[Remark~3.4]{Kat08b}. In contrast to the construction in \cite{Kat08b}, we obtain groupoids with totally disconnected unit spaces in this way.

For $i=0,1$, consider the exact sequence of $\Zz \Gamma$-modules
$$0 \to N_i \overset{\iota_i}{\lori} \Zz[K_i] \overset{\pi_i}{\lori} K_i \to 0,$$
where $\pi_i$ is the canonical projection and $N_i = \ker(\pi_i)$. All the modules in the sequence are countable. Since every Sylow subgroup of $\Gamma$ is cyclic, \cite[Theorem~1.4 and Corollary~3.7]{Kat08a} imply that $N_i$ is permutation projective for $i=0,1$. Hence, by \cite[Lemma~3.12]{Kat08a} and (the proof of) \cite[Lemma~3.2]{Kat08a}, we can find countable $\Gamma$-sets $X_i$ such that $N_i \oplus \Zz[X_i] \cong \Zz[X_i] \cong \Zz[K_i] \oplus \Zz[X_i]$ as $\Zz \Gamma$-modules. We obtain the exact sequence of $\Zz \Gamma$-modules
$$0 \to \Zz[X_i] \overset{\ti{\iota}_i}{\lori} \Zz[K_i \sqcup X_i] \overset{\ti{\pi}_i}{\lori} K_i \to 0,$$
where $\ti{\iota}_i$ is the composition $\Zz[X_i] \cong N_i \oplus \Zz[X_i] \overset{\iota_i \oplus \id}{\lori} \Zz[K_i] \oplus \Zz[X_i] \cong \Zz[K_i \sqcup X_i]$ and $\ti{\pi}_i$ is the composition $\Zz[K_i \sqcup X_i] \cong \Zz[K_i] \oplus \Zz[X_i] \overset{\pi_i + 0}{\lori} K_i$. Setting $A_i \defeq K_i \sqcup X_i$, $B_i \defeq X_i$, and $a_0 \defeq u \in K_0 \sqcup X_0 = A_0$, we obtain, for $i=0,1$, a short exact sequence of $\Zz \Gamma$-modules
$$0 \to \Zz[B_i] \to \Zz[A_i] \overset{\ti{\pi}_i}{\lori} K_i \to 0,$$
and $a_0 \in A_0$ is fixed by $\Gamma$ with $\ti{\pi}_0(a_0) = u$.

This is precisely the initial data needed to run the construction in \cite[Theorems~2.1 and 2.2]{Sp07b}. Let us sketch the construction and refer the reader to \cite{Sp07a,Sp07b} for details. Using \cite[Theorem~2.1]{Sp07b} and (the proof of) \cite[Theorem~2.2]{Sp07b}, we obtain irreducible (i.e., strongly connected) countable graphs $E_i$, $F_i$, $i=0,1$, each of which has a unique vertex emitting infinitely many edges, together with actions $\Gamma \curvearrowright E_i$, $i=0,1$, and $\Gamma$-equivariant inclusions $A_i \into E_i^0$, a vertex $w \in F_0^0$, such that we have isomorphisms
\begin{itemize}
\item $(K_0(\cO(E_i)),K_1(\cO(E_i)) \cong (K_i,\gekl{0})$, for $i=0,1$;
\item $(K_0(\cO(F_0)),[P_w],K_1(\cO(F_0))) \cong (\Zz,1,\gekl{0})$;
\item $(K_0(\cO(F_1)),K_1(\cO(F_1))) \cong (\gekl{0},\Zz)$;
\end{itemize}
and the first two isomorphisms fit into commutative diagrams, for $i=0,1$:
\[
\begin{tikzcd}
  A_i \arrow[d, "\ti{\pi}_i"] \arrow[r, hook] & E_i^0 \arrow[r] & \cP(\cO(E_i)) \arrow[d]
  \\
  K_i \arrow[rr, "\cong"] & & K_0(\cO(E_i))
\end{tikzcd}
\]
Note that commutativity of the diagram for $i=0$ implies that the isomorphism $K_0(\cO(E_0)) \cong K_0$ sends $[P_{a_0}]$ to $u \in K_0$. Now, the construction in \cite{Sp07a}, applied to $E_0$, $E_1$, $F_0$ and $F_1$, yields a \an{hybrid graph structure} $\Omega$. The $\Gamma$-actions on $E_0$ and $E_1$ induce a $\Gamma$-action on $\Omega$. In $\Omega$, we have a notion of finite paths, denoted by $\Omega^*$, and a length function $l: \: \Omega^* \to \Nz^2$ (where $\Nz = \gekl{0, 1, 2, \dotsc}$). We also have infinite paths in $\Omega$, and we can equip the set $Z$ of infinite paths with a natural topology such that $Z$ becomes a totally disconnected, second countable, locally compact Hausdorff space. Moreover, we can concatenate finite paths with infinite paths, written $\mu z$ for a finite path $\mu$ and an infinite path $z$. We can then define a groupoid $\ti{G}$ consisting of all triples $(z,n,z') \in Z \times \Zz^2 \times Z$ for which there exist finite paths $\mu$, $\nu$ and an infinite path $\zeta$ such that $z = \mu \zeta$, $z' = \nu \zeta$ and $l(\mu) - l(\nu) = n$. Inversion and multiplication are given by $(z,n,z')^{-1} = (z',-n,z)$ and $(z,m,z') (z',n,z'') = (z,m+n,z'')$. In particular, $\ti{G}^{(0)}$ is given by all triples $(z,0,z)$, $z \in Z$, and is canonically homeomorphic to $Z$. By construction, $\ti{G}$ is \etale{}, and it is topologically principal by \cite[Lemma~2.18]{Sp07a}. The $\Gamma$-action on $\Omega$ induces a $\Gamma$-action on $\ti{G}$. Spielberg shows in \cite{Sp07a} that $C^*_r(\ti{G})$ is a stable UCT Kirchberg algebra and that there is a canonical $\Gamma$-equivariant inclusion
$$
  i: \: (\cO(E_0) \otimes \cO(F_0)) \oplus (\cO(E_1) \otimes \cO(F_1)) \into C^*_r(\ti{G}),
$$
where $\cO(F_0)$ and $\cO(F_1)$ are equiped with the trivial $\Gamma$-actions. We do not need the precise form of $i$, but only the following properties:
\begin{enumerate}
\item[(i$_1$)] $i$ induces a $\Gamma$-equivariant isomorphism in K-theory;
\item[(i$_2$)] for every finite path $\mu \in E_0^*$ in $E_0$, $i(S_{\mu} S_{\mu}^* \otimes P_w)$ is of the form $1_Y$ for some compact open subspace $Y \subseteq Z$.
\end{enumerate}
In (i$_2$), $S_{\mu}$ is the partial isometry corresponding to $\mu$ in the graph C*-algebra $\cO(E_0)$, and $1_Y$ is the characteristic function of $Y$, viewed as an element in $C_0(Z) \cong C_0(\ti{G}^{(0)}) \subseteq C^*_r(\ti{G})$.

In particular, setting $1_X = i(P_{a_0} \otimes P_w)$, $i$ induces an isomorphism
$$(K_0(C^*_r(\ti{G})),[1_X],K_1(C^*_r(\ti{G}))) \cong (K_0,u,K_1).$$
As $\Gamma$ fixes $a_0$ and $w$, $X$ is fixed by $\Gamma$. Hence $\Gamma$ acts on $\ti{G}_X^X$, thus on $C^*_r(\ti{G}_X^X)$, and the canonical $\Gamma$-equivariant inclusion $C^*_r(\ti{G}_X^X) \cong 1_X C^*_r(\ti{G}) 1_X \into C^*_r(\ti{G})$ induces a $\Gamma$-equivariant isomorphism 
$$(K_0(C^*_r(\ti{G}_X^X)),[1],K_1(C^*_r(\ti{G}_X^X))) \cong (K_0,u,K_1).$$
Now set $G \defeq \ti{G}_X^X$. The unit space of $G$ is obviously given by $X$ and hence a totally disconnected compact space. $G$ is \etale{} and topologically principal as $\ti{G}$ has these properties.

Finally, let 
$$j: \: \cO(E_0) \otimes \cO(F_0) \into (\cO(E_0) \otimes \cO(F_0)) \oplus (\cO(E_1) \otimes \cO(F_1)) \overset{i}{\into} C^*_r(\ti{G})$$
be the restriction of $i$ to $\cO(E_0) \otimes \cO(F_0)$. Then $j_*$ is an isomorphism in $K_0$. It follows from the K-theory computation for graph C*-algebras and the K{\"u}nneth Formula (see \cite{Sch}) that $\menge{j_*[P_v \otimes P_w]}{v \in E_0^0}$ generates the $K_0$-group of $C^*_r(\ti{G})$. As $E_0$ is strongly connected and has a vertex emitting infinitely many edges, the image of $\cP(C(X)) \setminus \gekl{0} \to K_0(C^*_r(G))$ contains the group generated by $\menge{j_*[P_v \otimes P_w]}{v \in E_0^0}$ in $K_0(C^*_r(\ti{G})) \cong K_0(C^*_r(G))$. This shows that $\cP(C(X)) \setminus \gekl{0} \to K_0(C^*_r(G))$ is surjective. To see that $\cP(C(X)) \setminus \gekl{0,1} \to K_0(C^*_r(G))$ is surjective, it suffices to find a compact open subset $\emptyset \neq Y \subsetneq X$ such that $[1_Y] = [1_X]$. This however is easy using property (i$_2$), strong connectedness of $E_0$ and the existence of a vertex emitting infinitely many edges in $E_0$.
\eproof

By applying Lemma~\ref{Lem:Sp+Kat} and then stabilizing, we obtain the following:
\bcor
\label{Cor:Sp+Kat-stable}

Let $K_0$, $K_1$ be countable abelian groups. Let $\Gamma$ be a finite group such that every Sylow subgroup of $\Gamma$ is cyclic. Let $\Gamma \curvearrowright K_i$ be actions for $i=0,1$.

Then there exists an \etale{} topologically principal groupoid $G$ with totally disconnected unit space $X$ and an action $\Gamma \curvearrowright G$ such that
\begin{enumerate}
\item[(1$_{\rm st}$)] $C^*_r(G)$ is a stable UCT Kirchberg algebra.
\item[(2$_{\rm st}$)] We have a $\Gamma$-equivariant isomorphism $(K_0(C^*_r(G)),K_1(C^*_r(G))) \cong (K_0,K_1)$.
\item[(3$_{\rm st}$)] The canonical map $\cP(C_0(X)) \setminus \gekl{0} \to K_0(C^*_r(G))$ is surjective.
\end{enumerate}
\ecor

\bremark
A finite group $\Gamma$ has the property that every Sylow subgroup is cyclic if and only if $\Gamma$ admits a presentation
$$
\Gamma = \langle a,b ~|~ a^m = 1 = b^n,\ b^{-1}ab = a^r\rangle,
$$
where $m$ is odd, $0 \leq r < m$, $r^n \equiv 1 \mod m$ and $m$ and $n(r-1)$ are coprime, see for example \cite[10.1.10]{Ro}. In particular, every such group is an extension of one cyclic group by another.
\eremark

The induced actions $\Gamma \curvearrowright C^*_r(G)$ in Lemma~\ref{Lem:Sp+Kat} and Corollary~\ref{Cor:Sp+Kat-stable} preserve the canonical Cartan subalgebras $C(X)$ and $C_0(X)$, respectively. Therefore, Lemma~\ref{Spec(B) tot disc}, Lemma~\ref{Lem:Sp+Kat}, Corollary~\ref{Cor:Sp+Kat-stable} and the Kirchberg-Phillips classification theorem \cite{Kir,Phi} immediately imply
\blemma \label{lifts UCT Kirchberg}
Let $\Gamma$ be a finite group such that every Sylow subgroup is cyclic. Let $A$ be a UCT Kirchberg algebra and $\sigma:\Gamma \curvearrowright K_*(A)$ an action preserving $[1] \in K_0(A)$ if $A$ is unital. There exists an action $\alpha:\Gamma \curvearrowright A$ and an $\alpha$-invariant inverse semigroup $\cS \subseteq A$ of partial isometries that generates $A$ as a C*-algebra such that
\begin{itemize}
\item[(1)] $K_*(\alpha_\gamma) = \sigma_\gamma$ for all $\gamma \in \Gamma$;
\item[(2)] $(A,C^*(E))$ is a Cartan pair such that $C^*(E)$ has totally disconnected spectrum and $\alpha_\gamma(C^*(E)) = C^*(E)$ for all $\gamma \in \Gamma$;
\item[(3)] $\cP(C^*(E)) \setminus \gekl{0,1} \to K_0(A)$ is surjective if $A$ is unital;
\item[(4)] $\cP(C^*(E)) \setminus \gekl{0} \to K_0(A)$ is surjective if $A$ is stable.
\end{itemize}
\elemma

\bremark
\label{Rem3:tight}
If $A$ is unital in Lemma~\ref{lifts UCT Kirchberg}, then Remark~\ref{Rem1:tight} implies that $\cS \subseteq A$ can be chosen so that $\cS \into A$ induces an isomorphism $A \cong C^*_{tight}(\cS)$.
\eremark

We obtain the following result for actions with the Rokhlin property on unital UCT Kirchberg algebras, which follows immediately by Lemma~\ref{lifts UCT Kirchberg}, Remark~\ref{Rem3:tight} and Theorem~\ref{Cartan fixed general}.

\btheo
Let $\Gamma$ be a finite group such that every Sylow subgroup is cyclic. Let $A$ be a unital UCT Kirchberg algebra with $A \cong A \otimes M_{|\Gamma|^\infty}$. Let $\alpha:\Gamma \curvearrowright A$ be an action with the Rokhlin property.

Then there exists an $\alpha$-invariant inverse semigroup $\cS \subseteq A$ of partial isometries inducing an isomorphism $A \cong C^*_{tight}(\cS)$ such that $C^*(E)$ is a Cartan subalgebra with totally disconnected spectrum, $\alpha_\gamma(C^*(E)) = C^*(E)$ for all $\gamma \in \Gamma$, and $\cP(C^*(E)) \setminus \gekl{0,1} \to K_0(A)$ is surjective.
\etheo

\btheo \label{thm: approx repr}
Let $n \geq 2$. Let $A$ be a unital UCT Kirchberg algebra and $\alpha: \Zz_n \curvearrowright A$ be an outer approximately representable action. Assume that $A \rtimes_\alpha \Zz_n$ absorbs the UHF algebra $M_{n^\infty}$ tensorially. Then the following are equivalent:
\begin{itemize}
\item[(i)] $A \rtimes_\alpha \Zz_n$ satisfies the UCT;
\item[(ii)] there exists an inverse semigroup $\cS \subseteq A$ of partial isometries inducing an isomorphism $A \cong C^*_{tight}(\cS)$ such that 
\begin{itemize}
\item[(1)] $\cS$ is $\alpha$-homogeneous;
\item[(2)] $C^*(E)$ is a Cartan subalgebra with spectrum homeomorphic to the Cantor space in both $A^\alpha$ and $A$;
\item[(3)] unitaries in $A^\alpha$ witnessing approximate representability can be chosen in $C^*(E)$;
\end{itemize}
\item[(iii)] there exists a Cartan subalgebra $C \subseteq A$ such that $\alpha(C) = C$.
\end{itemize}
\etheo
\bproof
We prove that (i) implies (ii). As $\alpha$ is outer, $B := A \rtimes_\alpha \Zz_n$ is simple by \cite[Theorem~3.1]{Kis81} and purely infinite by \cite[Theorem~3]{J}. Denote by $\gamma:\Zz_n \curvearrowright B$ the dual action of $\alpha$, which has the Rokhlin property by Theorem~\ref{thm:duality Iz Naw}. Lemma~\ref{lifts UCT Kirchberg} together with Theorem~\ref{Cartan fixed general} and Remark~\ref{Rem2:tight} now yield an inverse semigroup $\tilde{\cS} \subseteq B \rtimes_\gamma \Zz_n$ of $\hat{\gamma}$-homogeneous partial isometries inducing an isomorphism $B \rtimes_\gamma \Zz_n \cong C^*_{tight}(\tilde{\cS})$ such that the sub-C*-algebra $C^*(E)$ generated by the corresponding idempotent semilattice is a Cartan subalgebra in both $B$ and $B \rtimes_\gamma \Zz_n$; in addition, unitaries in $(B \rtimes_\gamma \Zz_n)^{\hat{\gamma}} = B$ witnessing approximate representability of $\hat{\gamma}$ can be chosen in $C^*(E)$ and the map $\cP(C^*(E)) \setminus \gekl{0,1} \to K_0(B) = K_0((B \rtimes_\gamma \Zz_n)^{\hat{\gamma}})$ is surjective. Furthermore, being a masa with totally disconnected spectrum in a unital Kirchberg algebra, $C^*(E)$ has spectrum homeomorphic to the Cantor space. Thus, $\tilde{\cS} \subseteq B \rtimes_\gamma \Zz_n = A \rtimes_\alpha \Zz_n \rtimes_{\hat{\alpha}} \Zz_n$ satisfies the properties (1) - (3) for $\hat{\gamma}=\hat{\hat{\alpha}}$.

Now let $e \in \cB(\ell^{2}(\Zz_n))$ denote the rank one projection onto the subspace of constant functions and $\rho:\Zz_n \curvearrowright \cB(\ell^2(\Zz_n))$ the action induced by the right regular representation. There are equivariant isomorphisms
\[
(A,\alpha) \cong ((1 \otimes e)A \otimes \cB(\ell^2(\Zz_n))(1 \otimes e),\alpha \otimes \rho) \cong (f(A \rtimes_\alpha \Zz_n \rtimes_{\hat{\alpha}} \Zz_n)f,\hat{\hat{\alpha}})
\]
for some projection $f \in (A \rtimes_\alpha \Zz_n \rtimes_{\hat{\alpha}} \Zz_n)^{\hat{\hat{\alpha}}}\setminus\{0,1\}$, where the second isomorphism is induced by Takai duality, \cite[Theorem~3.4]{Tak}. Using that $\cP(C^*(E))\setminus \gekl{0,1} \to K_0(B) = K_0((A \rtimes_\alpha \Zz_n \rtimes_{\hat{\alpha}} \Zz_n)^{\hat{\hat{\alpha}}})$ is surjective,  we find a projection $p \in C^*(E) \setminus \gekl{0,1}$ such that $[p] = [f] \in K_0((A \rtimes_\alpha \Zz_n \rtimes_{\hat{\alpha}} \Zz_n)^{\hat{\hat{\alpha}}})$. Since $(A \rtimes_\alpha \Zz_n \rtimes_{\hat{\alpha}} \Zz_n)^{\hat{\hat{\alpha}}}$ is purely infinite simple, we find a unitary $u \in (A \rtimes_\alpha \Zz_n \rtimes_{\hat{\alpha}} \Zz_n)^{\hat{\hat{\alpha}}}$ such that $ufu^* = p$, see \cite[\S 1]{Cu1}. Hence 
\[
\Ad(u):(f (A \rtimes_\alpha \Zz_n \rtimes_{\hat{\alpha}} \Zz_n)f,\hat{\hat{\alpha}}) \stackrel{\cong}{\longrightarrow} (p(A \rtimes_\alpha \Zz_n \rtimes_{\hat{\alpha}} \Zz_n)p,\hat{\hat{\alpha}})
\]
is an equivariant isomorphism, showing that $(A,\alpha) \cong (p(A \rtimes_\alpha \Zz_n \rtimes_{\hat{\alpha}} \Zz_n)p,\hat{\hat{\alpha}})$. The inverse semigroup
\[
p \tilde{\cS} p \subseteq p(B \rtimes_\gamma \Zz_n)p = p(A \rtimes_\alpha \Zz_n \rtimes_{\hat{\alpha}} \Zz_n)p
\]
clearly consists of $\hat{\hat{\alpha}}$-homogeneous partial isometries, and $p \tilde{\cS} p \into p(A \rtimes_\alpha \Zz_n \rtimes_{\hat{\alpha}} \Zz_n)p$ induces an isomorphism $p(A \rtimes_\alpha \Zz_n \rtimes_{\hat{\alpha}} \Zz_n)p \cong C^*_{tight}(p \tilde{\cS} p)$ by Remark~\ref{Rem1:tight}. Moreover, $pC^*(E)p$ is a Cartan subalgebra in both  $p(A \rtimes_\alpha \Zz_n \rtimes_{\hat{\alpha}} \Zz_n)p$ and 
$$
p(A \rtimes_\alpha \Zz_n \rtimes_{\hat{\alpha}} \Zz_n)^{\hat{\hat{\alpha}}}p = (p(A \rtimes_\alpha \Zz_n \rtimes_{\hat{\alpha}} \Zz_n)p)^{\hat{\hat{\alpha}}},
$$
and coincides with the C*-algebra generated by the idempotent lattice of $p \tilde{\cS} p$. Furthermore, by multiplying unitaries in $C^*(E)$ witnessing approximate representability of $\hat{\hat{\alpha}}$ with $p$, one obtains unitaries in $pC^*(E)p$ witnessing approximate representability of the restricted action of $\hat{\hat{\alpha}}$ on $p(A \rtimes_\alpha \Zz_n \rtimes_{\hat{\alpha}} \Zz_n)p$. The desired inverse semigroup $\cS \subseteq A$ is now given as the image of $p\tilde{\cS}p$ under the equivariant isomorphism $(A,\alpha) \cong (p(A \rtimes_\alpha \Zz_n \rtimes_{\hat{\alpha}} \Zz_n)p,\hat{\hat{\alpha}})$.

The implication (ii) implies (iii) is trivial and the last implication, (iii) implies (i), follows from \cite[Proposition~3.4]{BL}.
\eproof

\bdefin [{\cite[Definition~3.6]{Iz}}]
Let $\Gamma$ be a finite abelian group and $A$ a unital, separable C*-algebra. An action $\alpha:\Gamma \curvearrowright A$ is called strongly approximately inner if there exist unitaries $u_n \in A^\alpha$, $n \in \Nz$ such that $\alpha = \lim_{n \to \infty} \Ad(u_n)$ in the topology of pointwise norm convergence.
\edefin

\bremark
$ $
\begin{itemize}
\item[1)]
By \cite[Proposition~4.8]{Sz2} and \cite[Corollary~4.10]{Sz1}, the assumption in Theorem~\ref{thm: approx repr} that $A \rtimes_\alpha \Zz_n$ absorbs $M_{n^\infty}$ tensorially is equivalent to $\alpha$ being conjugate to $\alpha \otimes \id:\Zz_n \curvearrowright A \otimes M_{n^\infty}$. On the other hand, if $\alpha:\Zz_n \curvearrowright A$ is any outer strongly approximately inner action on a unital UCT Kirchberg algebra, then $\alpha \otimes \id : \Zz_n \curvearrowright A \otimes M_{n^\infty}$ is outer and approximately representable by \cite[Lemma~3.10]{Iz}. Hence, the assumptions of Theorem~\ref{thm: approx repr} are satisfied for $\alpha \otimes \id : \Zz_n \curvearrowright A \otimes M_{n^\infty}$ as $(A \otimes M_{n^\infty})\rtimes_{\alpha \otimes \id} \Zz_n \cong (A \rtimes_\alpha \Zz_n) \otimes M_{n^\infty}$.
\item[2)] To the best of the authors' knowledge, it is not known whether all $\Zz_2$-actions on $\cO_2$ are strongly approximately inner. 
\end{itemize}
\eremark

We also record the following consequence for outer strongly approximately inner actions on $\cO_2$ by cyclic groups with prime power order.

\bcor \label{cor: charac UCT O_2}
Let $p \geq 2$ be a prime number, $n \geq 1$ a natural number and $q = p^n$. Let $\alpha: \Zz_q \curvearrowright \cO_2$ be an outer strongly approximately inner action. Then the following are equivalent 
\begin{itemize}
\item[(i)] $\cO_2 \rtimes_\alpha \Zz_q$ satisfies the UCT;
\item[(ii)] there exists an inverse semigroup $\cS \subseteq \cO_2$ of partial isometries inducing an isomorphism $\cO_2 \cong C^*_{tight}(\cS)$ such that 
\begin{itemize}
\item[(1)] $\cS$ is $\alpha$-homogeneous;
\item[(2)] $C^*(E)$ is a Cartan subalgebra with spectrum homeomorphic to the Cantor space in both $\cO_2^\alpha$ and $\cO_2$;
\item[(3)] unitaries in $\cO_2^\alpha$ witnessing approximate representability can be chosen in $C^*(E)$.
\end{itemize}
\item[(iii)] there exists some Cartan subalgebra $C \subseteq \cO_2$ such that $\alpha(C) = C$.
\end{itemize}
\ecor
\bproof
It follows from the Pimsner-Voiculescu sequence for $\Zz_q$-actions (see \cite[Theorem~10.7.1]{Bla}) that $K_*(\cO_2 \rtimes_\alpha \Zz_q)$ is uniquely $p$-divisible and in Cuntz standard form (see \cite[Lemma~4.4]{Iz} for details). Therefore $\cO_2 \rtimes_\alpha \Zz_q$ absorbs $M_{p^\infty} \cong M_{q^\infty}$ tensorially by the Kirchberg-Phillips classification theorem \cite{Kir, Phi}. Furthermore, $\alpha$ is approximately representable by \cite[Theorem~4.6]{Iz}. Hence, Theorem~\ref{thm: approx repr} applies and the proof is complete.
\eproof

The following result gives a characterization of the UCT problem for separable nuclear C*-algebras that are $KK$-equivalent to their $M_{p^\infty}$-stabilisation, where $p \geq 2$ is a prime number, in terms of outer $\Zz_p$-actions and Cartan subalgebras of $\cO_2$. A similar result for the case $p = 2$ can be found in \cite[Theorem~5.8]{BL}.

\btheo \label{thm: prime UCT problem}
Let $p \geq 2$ a prime number. The following statements are equivalent:
\begin{enumerate}
\item[(i)] Every separable nuclear C*-algebra $A$ that is $KK$-equivalent to $A \otimes M_{p^\infty}$ satisfies the UCT;
\item[(ii)] for every outer strongly approximately inner action $\alpha:\Zz_p \curvearrowright \cO_2$ there exists an inverse semigroup $\cS \subseteq \cO_2$ of partial isometries inducing an isomorphism $\cO_2 \cong C^*_{tight}(\cS)$ such that 
\begin{itemize}
\item[(1)] $\cS$ is $\alpha$-homogeneous;
\item[(2)] $C^*(E)$ is a Cartan subalgebra with spectrum homeomorphic to the Cantor space in both $\cO_2^\alpha$ and $\cO_2$;
\item[(3)] unitaries in $\cO_2^\alpha$ witnessing approximate representability can be chosen in $C^*(E)$;
\end{itemize}
\item[(iii)] every outer strongly approximately inner $\Zz_p$-action on $\cO_2$ fixes some Cartan subalgebra $B \subseteq \cO_2$ globally.
\end{enumerate}
\etheo
\bproof
Assume that (i) holds and let $\alpha: \Zz_p \curvearrowright \cO_2$ be an outer strongly approximately inner action. Then $\cO_2 \rtimes_\alpha \Zz_p$ is a unital, $M_{p^\infty}$-absorbing Kirchberg algebra. By assumption, $\cO_2 \rtimes_\alpha \Zz_p$ satisfies the UCT. Claim (ii) now follows from Corollary~\ref{cor: charac UCT O_2}.

The implication from (ii) to (iii) is trivial.

Lastly, assume (iii) and let $A$ be a separable, nuclear C*-algebra with the property that it is $KK$-equivalent to $A \otimes M_{p^\infty}$. Then $A$ is $KK$-equivalent to a unital $M_{p^{\infty}}$-absorbing Kirchberg algebra $A'$ by \cite[Theorem I]{Kir}. By \cite[Proposition~4.14]{BS}, there exists an outer approximately representable action $\gamma: \Zz_p \curvearrowright \cO_2$ such that $\cO_2 \rtimes_\gamma \Zz_p$ is $KK$-equivalent to $M_{p^{\infty}}^{\oplus p-1}$. Fix an isomorphism $A' \otimes \cO_2 \cong \cO_2$ and let $\alpha$ be the automorphism corresponding to $\id_{A'} \otimes \gamma$ under this identification. Then
\[
\cO_2 \rtimes_\alpha \Zz_p \cong (A' \otimes \cO_2) \rtimes_{\id_{A'} \otimes \gamma} \Zz_p \sim_{KK} A' \otimes M_{p^\infty}^{\oplus p-1} \cong A'^{\oplus p-1} \sim_{KK} A^{\oplus p-1}.
\]
As an outer strongly approximately inner action on $\cO_2$, $\alpha$ fixes some Cartan subalgebra $B \subseteq \cO_2$ globally by assumption. By \cite[Proposition~3.4]{BL}, $\cO_2 \rtimes_\alpha \Zz_p$ satisfies the UCT and so does $A$. This shows (i) and the proof is complete.
\eproof

The UCT problem can now be reformulated as follows.

\bcor
\label{Cor:UCTEquivChar}
The following statements are equivalent:
\begin{enumerate}
\item[(i)] Every separable nuclear C*-algebra $A$ satisfies the UCT;
\item[(ii)] for every prime number $p \geq 2$ and every outer strongly approximately inner action $\alpha:\Zz_p \curvearrowright \cO_2$ there exists an inverse semigroup $\cS \subseteq \cO_2$ of partial isometries inducing an isomorphism $\cO_2 \cong C^*_{tight}(\cS)$ such that
\begin{itemize}
\item[(1)] $\cS$ is $\alpha$-homogeneous;
\item[(2)] $C^*(E)$ is a Cartan subalgebra with spectrum homeomorphic to the Cantor space in both $\cO_2^\alpha$ and $\cO_2$;
\item[(3)] unitaries in $\cO_2^\alpha$ witnessing approximate representability can be chosen in $C^*(E)$;
\end{itemize}
\item[(iii)] every outer strongly approximately inner $\Zz_p$-action on $\cO_2$ with $p=2$ or $p = 3$ fixes some Cartan subalgebra $B \subseteq \cO_2$ globally.
\end{enumerate}
\ecor
\bproof
This follows immediately from Theorem~\ref{thm: prime UCT problem} and \cite[Proposition~4.16]{BS}.
\eproof

\bremark
The reader may find a comment (without further details) in \cite[23.15.12]{Bla} pointing towards similar characterizations of the UCT problem as those in Theorem~\ref{thm: prime UCT problem} and Corollary~\ref{Cor:UCTEquivChar}.
\eremark

\end{document}